\documentclass[12pt]{amsart}

\usepackage[left=2.7cm,right=2.7cm,top=2.9cm,bottom=2.9cm]{geometry}

\usepackage[utf8]{inputenc}
\usepackage{graphicx,color}
\usepackage{latexsym}
\usepackage{graphicx} 
\usepackage{esint} 
\usepackage{amsthm,amsmath, amssymb,amsfonts,esint} 
\usepackage{layout,textcomp}
\usepackage{enumerate}

\newenvironment{acknowledgement}{\noindent\textbf{Acknowledgments: }}

\usepackage[colorlinks=true,urlcolor=blue,
citecolor=red,linkcolor=blue,linktocpage,pdfpagelabels,
bookmarksnumbered,bookmarksopen]{hyperref}

\usepackage[english]{babel}

\theoremstyle{plain}

\newcommand{\R}{\mathbb{R}}

\newtheorem{theorem}{Theorem}[section]

\newtheorem{lemma}[theorem]{Lemma}

\newtheorem{proposition}[theorem]{Proposition}

\theoremstyle{remark}

\numberwithin{equation}{section}

\title[Blow-up analysis of Large conformal metrics]{Blow-up analysis of Large conformal metrics with prescribed  Gaussian and geodesic curvatures}

\author[R. Caju]{Rayssa Caju}
\address[R. Caju]{Department of Mathematical Engineering, University of Chile
\newline\indent 
    Beauchef 851, Edificio Norte, Santiago, Chile}
\email{\href{mailto:rayssacaju@gmail.com}{rcaju@dim.uchile.cl}}

\author[T. Cruz]{Tiarlos Cruz}
\address[T. Cruz]{ Institute of Mathematics, 
	Federal University of Alagoas
	\newline\indent 
	57072-970, Maceió-AL, Brazil}
\email{\href{mailto: cicero.cruz@im.ufal.br}{cicero.cruz@im.ufal.br}}

\author[A. Silva Santos]{Almir Silva Santos}

\address[A. Silva Santos]{Department of Mathematics, 
	Federal University of Sergipe
	\newline\indent 
	49107-230, Sao Cristóv\~ao-SE, Brazil}
\email{\href{mailto: almir@mat.ufs.br}{ almir@mat.ufs.br}}

\thanks{RC was partially supported by FONDECYT grant number 11230872 and by Centro de Modelamiento Matemático (CMM) BASAL fund FB210005 for center of excellence from ANID-Chile. TC was partially supported by CNPq grant number  307419/2022-3 and 405468/2021-0. ASS was partially supported by CNPq grant number 403349/2021-4, 408834/2023-4 and 312027/2023-0. The authors were  supported by FAPEAL, Brazil (Process E:60030.0000002329/202) and CNPq grant number 408834/2023-4
}

\begin{document}
\begin{abstract}
Consider a compact Riemannian surface $(M,g)$  with a nonempty boundary and negative Euler characteristic. Given two smooth non-constant functions $f$ in $M$ and $h$ in $\partial M$ with $\max f= \max h= 0$, under a suitable condition on the maximum points of $f$ and $h$, we prove that for sufficiently small positive constants $\lambda$ and $\mu$, there exist at least two distinct conformal metrics $g_{\lambda,\mu}=e^{2u_{\mu,\lambda}}g$ and $g^{\lambda,\mu}=e^{2u^{\mu,\lambda}}g$ with prescribed sign-changing  Gaussian and geodesic curvature equal to $f + \mu$ and  $h + \lambda,$ respectively. Additionally, we employ the method used by Borer et al. (2015) to study the blowing-up behavior of the large solution $u^{\mu,\lambda}$ when $\mu\downarrow 0$ and $\lambda\downarrow 0$. Finally, we derive a new Liouville-type result for the half-space, eliminating one of the potential blow-up profiles.
\end{abstract}

\maketitle

\section{Introduction}

The problem of prescribing curvature has attracted much attention in the last few decades due to its relevance and importance in various branches of mathematics.
Up to our knowledge, M. Berger \cite{MR0295261} was the first to study the problem of prescribing the Gaussian curvature on a given closed surface with negative Euler characteristic. Later, in a series of papers \cite{MR0675736,MR0343205,MR0375154,MR0375153,MR0365409}, J. Kazdan and F. Warner have addressed the problem of prescribing the Gaussian and scalar curvature on a closed manifold, in the conformal and nonconformal setting. It turned out that conformal deformation of the metric is an important tool for finding geometric properties by using analysis techniques.

Let $(M,g)$ be a closed Riemannian surface. Through a conformal deformation of the metric $\overline g=e^{2u}g,$ for some smooth function $u$, it is well known that the Gaussian curvature $K_g$ and $K_{\overline g}$ 
of $g$ and $\overline g$, respectively, are related by the following PDE
\begin{equation}\label{eq053}
    -\Delta_gu + K_{g} = K_{\bar{g}}e^{2u},
\end{equation}
where $\Delta_g$ denotes the Laplace-Beltrami operator with respect to $g.$
According to the Gauss-Bonnet Theorem, the solvability of equation \eqref{eq053} depends on the sign of the Euler characteristic $\chi(M)$, which imposes necessary conditions on the function to be realized as the Gaussian curvature of $\overline g$. A complete understanding of which functions can be the Gaussian curvature of some conformal metric remains a challenging problem, especially in the case of the two-sphere, commonly referred to as the Nirenberg problem. We direct the reader to the works \cite{MR4277273, MR3694645} and the references therein.

Consider a compact Riemannian surface $(M,g)$ with nonempty boundary $\partial M$. If $\overline g$ is a conformal metric to $g$ with Gaussian and geodesic curvature equal to $f\in C^\infty(M)$ and $h\in C^\infty(\partial M)$, respectively, then there exists a smooth function $u$ such that $\overline g=e^{2u}g,$ which satisfies the boundary value problem
\begin{equation}\label{eq000}
\begin{cases}
    -\Delta_g u+K_g=fe^{2u} & \mbox{ in }M,\\
    \dfrac{\partial u}{\partial\nu_g}+\kappa_g=he^u & \mbox{ on }\partial M,
\end{cases}
\end{equation}
where $\nu_g$ is the outward unit normal to  $\partial M$  and $\kappa_g$ is the geodesic curvatures of the boundary with respect to  $g$.

Although the existence of solutions to \eqref{eq000} is not completely understood, there are many works dedicated to studying the existence and also the blow-up analysis of these solutions. Since the literature on this topic is extensive, we only refer to recent works  \cite{MR4517687, lopez2024prescribing,struwe2024prescribed} and the references therein.

If \eqref{eq000} has a solution, by the Gauss-Bonnet Theorem, we obtain
$$\int_{M}fdv_{\overline g} + \int_{\partial M}hd\sigma_{\overline g} = 2\pi\chi(M), $$
where $dv_{\overline g}=e^{2u}dv_g$ and $d\sigma_{\overline g}=e^ud\sigma_g$ are the area and length elements, respectively, of the metric $\overline g=e^{2u}g$. This means that the solubility of \eqref{eq000} is also subject to topological constraints. For example, if $\chi(M)<0$, $f$ and $h$ are non-negative functions, then there is no solution.

The main goal of this work is to study the structure of the set of solutions of \eqref{eq000} in the case $\chi(M)<0$. We will extend the results for closed surfaces obtained in \cite{MR3351750} to compact surfaces with nonempty boundary, as well as the corresponding improvements provided in \cite{MR4153107} of the asymptotic behavior of the solutions obtained in \cite{MR3351750}. In future work, we will address the case $\chi(M)=0$.

Without loss of generality,  by the celebrated uniformization theorem of B. Osgood, R. Phillips, and P. Sarnak for surfaces with boundary \cite{MR0960228} we may assume that $K_g=-1$ and $\kappa_g\equiv 0$ (see also  \cite[Proposition 3.1]{MR4517687} and  \cite[Lemma 2.1]{MR4225836}). Taking this into account, solving  equation \eqref{eq000} is equivalent to solving 
\begin{equation}\label{eq001}
    \begin{cases}
    -\Delta_gu-1=fe^{2u} & \mbox{ in }M,\\
    \dfrac{\partial u}{\partial\nu_g}=he^u & \mbox{ on }\partial M.
\end{cases}
\end{equation}

The problem \eqref{eq001} is the Euler-Lagrange equation of the functional $I_{f,h}:H^1(M)\rightarrow\mathbb R$ given by
\begin{equation}\label{functional}I_{f,h}(u)=\frac{1}{2}\int_{M}\left(|\nabla u|^2-2u-fe^{2u}\right)dv_g-\int_{\partial M}he^u d\sigma_g.
\end{equation}
Finding a critical point to $I_{f,h}$, for general functions $f$ and $h$, is a challenging task. The area and boundary term compete, and it may not be possible to verify whether $I_{f,h}$ is bounded from below or not. For non-positive functions $f$ and $h$ with one of them being non-identically zero, \cite[Theorem 2]{MR0749522} implies the existence of a solution to \eqref{eq001}. Moreover, in this case, it is not difficult to see that the functional \eqref{functional} is strictly convex (see also \cite[Theorem 1.1]{MR4517687}). Thus, it holds the following result.
\begin{theorem}
Let $(M,g)$ be a compact Riemannian surface with $\partial M\not=\emptyset$ and $\mathcal X(M)<0$. If $f\in C^\infty(M)$ and $h\in C^\infty(\partial M)$ are nonpositive functions with $f\not\equiv 0$ or $h\not\equiv 0$, then \eqref{eq000} admits a unique solution.
\end{theorem}

Recently, existence of solutions to \eqref{eq001} was studied by R. López-Soriano, A. Malchiodi, and D. Ruiz \cite{MR4517687}. They proved the existence of solutions when $\mathcal X(M)<0$ under the assumptions $f<0$ and  $h(p)/\sqrt{|f(p)|}<1$ for all $p\in\partial M$. Moreover, if $h\leq 0$, then the solution is unique, see \cite[Theorem 1.1]{MR4517687}. Their proof consists of showing  that the functional is coercive, weakly lower semicontinuous, and strictly convex when $h\leq 0$. They also studied existence results of solution for \eqref{eq000} in the case $\mathcal X(M)=0$, under extra assumptions on the integral of $h$, and analyzed the blow-up phenomenon. See \cite{lopez2024prescribing} for results in the disk case.

We are interested in studying the solution set of \eqref{eq001} when  $f$ and $h$ change signs. First, considering the variational approach, since we are interested in critical points for $I_{f,h}$, it is natural to study nondegeneracy of its local minimizer, which is determined by the second derivative
\begin{equation}\label{eq047}
d^2I_{f,h}(u)(m,m) = \int_M(|\nabla m|^2 - 2fm^2e^{2u})dv_g - \int_{\partial M}hm^2e^{u}d\sigma_g.
\end{equation}

Building upon the ideas in  \cite{MR3628322} (see also \cite{MR1753266}), we can prove our first result.

\begin{theorem}\label{Theo_exist}
Let $(M,g)$ be a compact Riemannian surface with $\partial M\not=\emptyset$ and $\chi(M)<0$.  Suppose that for $f\in C^{0,\alpha}(M)$ and $h\in C^{1,\alpha}(\partial M)$ with $h\leq 0$, for some $0<\alpha<1,$ the functional $I_{f,h}$  admits a local minimizer $u_{f,h}\in H^1(M)$. Then
\begin{enumerate}[(a)]
    \item there exists a constant $c_0>0$ such that 
    $$d^2I_{f,h}(u_{f,h})(m,m)\geq c_0\|m\|^2_{H^1},$$
    for all  $m\in H^1(M)$. In particular,  $u_{f,h}$ is a non-degenerate critical point of $I_{f,h}$. 

 \item there exists a neighborhood $U\times W\subset C^{0,\alpha}(M)\times C^{1,\alpha}(\partial M)$ of $(f,h)$ such that for any  $(z,k)\in U\times W$,  there exists a strict local minimizer for  $I_{z,k}$ in $C^{2,\alpha}(M)$  smoothly dependent on $z,k$. In particular, if $f$ and  $h$  are smooth, then the minimizer is smooth as well.
\end{enumerate}
\end{theorem}

 Item (a) can be viewed as the counterpart to some results in \cite{MR1305205,MR3351750} in the context of surfaces with boundary. It would be interesting to know whether this result holds for general $h\in C^\infty(\partial M)$. Item (b) follows as a consequence of the implicit function theorem.

 We now describe an interesting phenomenon. Consider nonconstant functions $f\in C^\infty(M)$ and $h\in C^\infty(\partial M)$ with $\max f=\max h=0$. For real numbers $\mu>0$ and $\lambda>0,$ define  $f_{\mu}=f+\mu$ and $h_{\lambda} = h+\lambda$. First, we observe that if  $\mu$ and $\lambda$ are large enough, it follows from the Gauss-Bonnet Theorem that \eqref{eq001} has no solution. On the other hand, by Theorem \ref{Theo_exist}, for $\mu$ and $\lambda$ small enough, the equation \eqref{eq001} admits a smooth solution $u_{\mu,\lambda}$, which is a strictly local minimizer of  $I_{f_\mu,h_\lambda}$. But since $\max f_\mu>0$, it is not difficult to see that the functional $I_{f_\mu,h_\lambda}$ is no longer bounded from below, and we may expect the existence of a further critical point $u^{\mu,\lambda}$ of saddle-type.  In this regard, we obtain the following result.

\begin{theorem}\label{second_critical}
    Let $(M,g)$ be a compact Riemannian surface with $\partial M\not=\emptyset$ and $\chi(M)<0$. Let $f\in C^\infty(M)$ and $h\in C^\infty(\partial M)$ be nonconstant functions satisfying $\max f=\max h=0$. Consider the families of functions  $h_\lambda=h+\lambda$  and $f_\mu=f+\mu,$  where $\lambda$ and $\mu$ are two real numbers. There exist real numbers $\lambda_{0}$ and $\mu_0$ such that for almost every $0<\lambda<\lambda_{0}$ and  $0<\mu<\mu_{0}$, the functional $I_{f_\mu,h_\lambda}$ admits a strict local minimizer $u_{\lambda,\mu}$ and a second  critical point $u^{\mu,\lambda}\neq u_{\mu,\lambda}$ of mountain-pass type.

\end{theorem}

    In the context of closed surfaces, the existence of a second solution to \eqref{eq053}  was first proved in \cite{MR1257102}, when $\chi (M)<0$ and the prescribed Gaussian curvature $f_\lambda=f+\lambda$, with $\lambda>0$ small enough, changes sign. However, the method used does not provide the geometric shape of the solutions. Later, in \cite{MR3351750}, the authors gave a new proof using the so-called ``monotonicity trick" (see \cite{MR0970849, MR0926524} and also \cite{MR1619043}), which allowed them to bound the volume of the ``large" solutions when $\lambda\to 0$ and then perform a blow-up analysis. The blow-up limit is spherical, as proved in a further work by M. Struwe \cite{MR4153107}. See also \cite{MR3385180}, where, via a refined analysis, the authors constructed ``branches" of metrics exhibiting multiple spherical bubbling in the limit solution by matched asymptotical expansion. In \cite{MR3688436} the authors studied large conformal metrics with prescribed scalar curvature in a similar manner as in \cite{MR3351750,MR3385180}.

In contrast to what happens in the closed surface case, the presence of a boundary crucially gives rise to a new non-compactness phenomenon, especially the existence of blow-up profiles with infinite area, as shown in \cite{MR4517687}. In this work, the authors exploit the variational formulation of the problem to obtain existence results, employing a combination of minimization and min-max techniques. We refer to \cite{MR3198871,battaglia2023mean,MR4591853,MR4553958,MR3399138,MR2274940,MR3249809,lopez2024prescribing,MR3892074,MR3767960}   for more references on blow-up analysis of solutions
and to \cite{MR4546499, MR4517687}  and the references therein for the most recent progress in the topic and for including new general blow-up results.

Our next goal is to understand what happens with the solution $u^{\mu,\lambda}$, given by Theorem \ref{second_critical}, when $\lambda\downarrow 0$ and $\mu\downarrow 0$. The main result in this respect is the following.

\begin{theorem}\label{theor-blow-up-analysis}
Let $(M,g)$ be a compact Riemannian surface with $\partial M\not=\emptyset$ and $\chi(M)<0$. Let $f\in C^\infty(M)$ and $h\in C^\infty(\partial M)$ be nonconstant functions with $\max f=\max h=0$. Suppose the maximum points of $f$ and $h$ are the same, and non-degenerate, and critical points of $f$. For $(\mu,\lambda)\in\mathbb R^2,$ define $f_\mu:=f+\mu$ and $h_\lambda:=h+\lambda$. Then, there exist two sequences of real numbers, $\mu_n\downarrow 0$ and $\lambda_n\downarrow 0$, 
and a sequence of solutions $u_n$ of 
\begin{equation}\label{eq017}
    \begin{cases}
    -\Delta_gu_n+K_g=f_{\mu_n} e^{2u_n} & \mbox{ in }M,\\
    \dfrac{\partial u_n}{\partial\nu_g}=h_{\lambda_n} e^{u_n} & \mbox{ on }\partial M,
\end{cases}
\end{equation}
which are non-minimizing critical points of $I_{f_{\mu_n},h_{\lambda_n}}$, such that the conformal metric  $g_n = e^{2u_n}g$ has total curvature uniformly bounded,
$$
\int_M|K_{g_n}|dv_{g_n}+\int_{\partial M}|\kappa_{g_n}|d\sigma_{g_n} <C <\infty.
$$
Here $K_{g_n}=f_{\mu_n}$ and $\kappa_{g_n}=h_{\lambda_n}$ are the Gaussian and geodesic curvature of the metric $g_n$, respectively. Moreover, there exist at most three points $\{p_\infty^{(i)}\in M:i=1,\ldots,N\}$ with $N\in\{1,2,3\}$ and $f(p_\infty^{(i)})=h(p_\infty^{(i)})=0$, and two sequences $r_n^{(i)}\downarrow 0$ and $p_n^{(i)}\rightarrow p_\infty^{(i)}$ in $M$, such that
\begin{enumerate}[(a)]
\item $u_n\rightarrow u_\infty$ smooth locally on $M_\infty:=M\backslash\{p_\infty^{(i)}:1\leq i\leq N\}$, and $u_\infty$ induces a complete metric $g_\infty:=e^{2u_\infty}g$ on $M_\infty$ of finite total curvature $$\int_{M_\infty}|K_{g_\infty}|dv_{g_\infty}+\int_{\partial M_\infty}|\kappa_\infty|d\sigma_{g_\infty}<\infty,$$
where $K_{g_\infty}=f$ and $\kappa_{g_\infty}=h$ are the Gaussian and geodesic curvatures of $g_\infty$, respectively.
    \item For each $i\in\{1\ldots, N\},$ there holds $r_n^{(i)}/\lambda_n\rightarrow 0$ and in local conformal coordinates $\varphi_i:U_i\subset\mathbb R^2\to M$ around $p_\infty^{(i)}$, with $\varphi_i^*g=e^{2v^{(i)}}g_{euc}$, where $g_{euc}$ is the euclidean metric in $\mathbb R^2$, we have that
        $$w_{n}(x):= v_{n}(x_n+r_{n}^{(i)}x) - v_{n}(x_n)\to w_\infty(x)=\log \frac{2\Lambda}{c_\infty\Lambda^2+\left(s-s_0\right)^2+\left(t-t_0+d_\infty\Lambda\right)^2},$$
        smoothly locally in $\mathbb R^2_{t_0}=\mathbb R\times (t_0,+\infty)$, for some $t_0\geq 0$, where $x=(s,t)$. Here,  $x_n=\varphi_i^{-1}(p_n^{(i)})$ and $v_n=u_n\circ\varphi_i+v^{(i)}$. In addition, $s_0\in\mathbb R$,  $\Lambda>0$ and $c_\infty,d_\infty\in [0,1]$, with $d_\infty=1$ if $c_\infty=0.$ The function $w_\infty$ induces a metric $g_\infty=e^{2 w_\infty}g_{euc}$ of Gaussian curvature $K_{g_\infty}=c_\infty$ on $\mathbb R^2_{t_0}$ and geodesic curvature $\kappa_{g_\infty}=d_\infty$ on $\partial \mathbb R^2_{t_0}$.

        Moreover, if $c_\infty=0$, then $N=1$. If $c_\infty\in (0,1]$ and $d_\infty=0$, then $1\leq N\leq 2$. If $c_\infty, d_\infty\in (0,1],$ then $1\leq N\leq 3$.
\end{enumerate}

\end{theorem}

Although the blow-up analysis follows a similar methodology to that of the case of closed surfaces, the presence of the boundary fundamentally changes the dynamics of our analysis. Primarily, there is a deep interdependence between the parameters $\mu$ and $\lambda$, as one needs to be correlated almost quadratically with the other. Secondly, differently from the closed case, by analyzing a volume estimate, we manage to exclude the possibility of whole bubbles of large curvature being formed in blow-up points happening in the boundary. Finally, following the ideas developed in the paper of \cite{MR4517687}, which majorly uses Ekeland’s variational principle, we are able to use the techniques applied in \cite{MR3351750,MR4753069,MR4153107} to understand the possibility of slow blow-ups. For technical reasons, it was necessary to assume in Theorem \ref{theor-blow-up-analysis}
that the functions $f$ and $h$ have the same maximum points. It would be nice to know what would happen if this is not the case.

We would like to emphasize that, although there are characterizations of blow-up solutions to the Neumann boundary value problem in the literature, as mentioned earlier, there are some issues in applying these results to prove Theorem \ref{theor-blow-up-analysis}. The first issue is that to use many of these works (see again \cite{MR3198871,MR4753069,MR3892074,MR3767960}), it is necessary to require that the functions which are expected to prescribe as the Gaussian and geodesic curvatures do not change sign, a condition that clearly does not hold in our case. The second is that these works do not allow to uniformly bound  the area and the boundary length of the ``large” solutions  $u_{\mu,\lambda}$ as $\mu,\lambda\to 0$ appropriately.  Fortunately, we overcome these obstacles by adapting the mountain pass technique used by Struwe \cite{MR0926524,MR1245097} in a way similar to \cite{MR3351750,MR4153107}. See \cite{MR3412381} for a similar result for the torus context, and \cite{MR3628322} for the $Q$-curvature context.

To the best of our knowledge, there are few results in the literature on large conformal metrics that simultaneously address the Gaussian and geodesic curvatures. Recently, a result in the disk has appeared, see  \cite{MR4204565}. In this work, the authors proved the existence of bubbling type solutions under certain nondegeneracy assumptions on their curvatures using a Lyapunov-Schmidt reduction.

The analysis done in \cite{MR3351750} gives two possibilities for the blow-up: one for the bubble solution, and another one for a solution that gives rise to a metric $g_\infty$ in all $\mathbb R^2$ with finite volume, finite total curvature and Gaussian curvature $K_{g_\infty}=1+(Ax,x)$, where $A$ is a negative definite $2\times 2$ matrix. Later, M. Struwe \cite{MR4153107} ruled out the second possibility by showing a Liouville-type theorem in the plane. See \cite{arXiv:2401.03457,MR4753069,MR4287915} for an extension of this theorem to higher-order operator. The proof of Theorem \ref{theor-blow-up-analysis} follows the same line. To rule out a  potential second scenario for the blow-up, we prove the following  Liouville-type theorem in the half-space $\mathbb R^2_+$, which might be of independent interest.

\begin{theorem}\label{teo001}
Suppose that $A$ is a negative definite $2\times 2$ matrix.  Then there is no solution $w \in C^{\infty}(\mathbb{R}^{2}_{+})$ of the equation
\begin{equation*}
\begin{cases}
    -\Delta w=Fe^{2w} & \mbox{ in }\R^{2}_{+},\\
    \dfrac{\partial w}{\partial\nu}=e^w & \mbox{ on }\partial \mathbb R^{2}_{+},
\end{cases}
\end{equation*}
with $w\leq C$, where $F(x):=1+ (Ax,x)$, such that $Fe^{2w}\in L^1(\mathbb R^2_+)$,
\begin{equation*}
 \int_{\R^{2}_{+}}e^{2w} dx<\infty,\quad  \int_{\partial\R^{2}_{+}}e^{w} dl<\infty,\quad \mbox{ and }\quad    \int_{\mathbb{R}^{2}_{+}}Fe^{2w} dx\in \R.
\end{equation*}
\end{theorem}

The paper is structured as follows. In Section \ref{sec:Stability-results}, we investigate the stability and nondegeneracy properties of local minimizers  with respect to the energy functional  \eqref{functional}. Next, in Section \ref{sec:existence-of-a-further}, we establish an existence result of a class of large conformal solutions of saddle-type as described in Theorem \ref{second_critical}. Moving forward to Section \ref{sec:blow-up-analysis}, in order to prove Theorem \ref{theor-blow-up-analysis}, we conduct a blow-up analysis to study the behavior of these large solutions as the parameters approach zero suitable. In Section \ref{rule} we prove Theorem \ref{teo001}, which finalizes the proof of Theorem \ref{theor-blow-up-analysis}. 

\medskip

\acknowledgement{This work was carried out while TC was visiting ICTP and ASS was visiting University of Chile.  ASS is very grateful for the hospitality of University of Chile.
TC is very grateful to Claudio Arezzo for the hospitality and also he would like to acknowledge support from the ICTP through the Associates Program (2018-2023). The authors would like to thank S. Almaraz and R. López-Soriano for useful comments on an earlier version of this paper. Finally, we thank the anonymous referee for its invaluable suggestions and corrections.}

\section{Stability  results}\label{sec:Stability-results}
Throughout this paper, $(M,g)$ will denote a compact, connected Riemannian surface with nonempty boundary and negative Euler characteristic, i.e., $\chi(M)<0$. In this case, the Gauss-Bonnet formula leads to
$$\int_M K_g  dv_g + \int_{\partial M} \kappa_g  d\sigma_g = 2\pi \chi(M) < 0,$$
where $K_g$ and $\kappa_g$ are the Gaussian and geodesic curvature, respectively, $dv_g$ and $d\sigma_g$ are the area and length element of $g$.

\subsection{Nondegeneracy of local minimizers}
We will guarantee the nondegeneracy of local minimizers for the energy $I_{f,h}$ for any arbitrary functions $f$ and $h$ such that $h \leq 0$, in the case that $\chi(M) < 0$. The next proposition proves item (a) of Theorem \ref{Theo_exist}.

\begin{proposition}\label{propo001}
Let $(M,g)$ be a compact Riemannian surface with $\partial M\not=\emptyset$ and $\chi(M)<0$. Let $f\in C^{0,\alpha}(M)$ and $h\in C^{0,\alpha}(\partial M)$ with $h\leq 0$. Suppose the functional $I_{f,h}$ given by \eqref{functional} admits a local minimizer $u_{f,h}\in H^1(M)$. Then $u_{f,h}$ is a non-degenerate critical point of $I_{f,h}$ in  the sense that $d^2I_{f,h}(u_{f,h})(m,m)\geq c_0\|m\|^2_{H^1},$
    for some positive constant $c_0$ and for all $m\in H^1(M)$.
\end{proposition}
\begin{proof}
Define $c_0:=\inf\{d^2I_{f,h}(u_{f,h})(m,m):\|m\|_{H^1}=1\}$. Since $u_{f,h} \in H^1(M)$ is a local minimizer of $I_{f,h}$, it follows that  $d^2I_{f,h}(u_{f,h})(m,m)\geq 0$ for all $m\in H^1(M)$. Thus $c_0 \geq 0$. To establish the desired result, we need to prove that $c_0 > 0$. By contradiction, let us assume that $c_0 = 0$.

\medskip

\noindent{\bf Claim 1:} There is $m\in H^1(M)$ with $\|m\|_{H^1}=1$ such that $d^2I_{f,h}(u_{f,h})(m,m)=0.$

\medskip

    Let $\{m_k\}_{k\in\mathbb N}$ be a sequence of functions in $H^1(M)$ with $\|m_k\|_{H^1}=1$ and such that $d^2I_{f,h}(u_{f,h})(m_k,m_k)\rightarrow 0$ as $k\rightarrow\infty$. Since a ball in any Hilbert space is weakly compact, there exists a function $m\in H^1(M)$ and  a subsequence of $\{m_k\}$, still denoted by $\{m_k\}$, converging weakly to $m$ in $H^1(M)$ and strongly in $L^p(M)$ for any $1\leq p<\infty$. This implies that 
    \begin{equation}\label{eq002}
    \|m\|_{H^1}\leq \liminf\|m_k\|_{H^1}.    
    \end{equation}

By a Moser-Trudinger inequality, see \cite{MR3836128,MR2178789}, and the trace embedding \cite[Theorem 6.2]{MR3014461}, see also \cite[Section 3.2]{phdthesis}, we obtain that $e^{u_{f,h}}$ belongs to $L^p(M)$ and $L^p(\partial M)$ for all $1\leq p<\infty$. From this and H\"older inequality, we get that  $fm_k^2e^{2u_{f,h}}\rightarrow fm^2e^{2u_{f,h}}$ in $L^1(M)$ and $hm_k^2e^{u_{f,h}}\rightarrow hm^2e^{u_{f,h}}$ in $L^1(\partial M)$. Using \eqref{eq047} and the fact that $d^2I_{f,h}(u_{f,h})(m,m)\geq 0$ we get
    \begin{align*}
        \|\nabla m_k\|^2_{L^2} & =d^2I_{f,h}(u_{f,h})(m_k,m_k)+2\int_Mfm_k^2e^{2u_{f,h}}dv_g+\int_{\partial M}hm_k^2e^{u_{f,h}}d\sigma_g\\
        & \rightarrow 2\int_Mfm^2e^{2u_{f,h}}dv_g+\int_{\partial M}hm^2e^{u_{f,h}}d\sigma_g\\
        & \leq d^2I_{f,h}(u_{f,h})(m,m)+2\int_Mfm^2e^{2u_{f,h}}dv_g+\int_{\partial M}hm^2e^{u_{f,h}}d\sigma_g= \|\nabla m\|^2_{L^2}.
    \end{align*}
    This implies that
    $$\limsup\|m_k\|^2_{H^1}=\limsup(\|m_k\|^2_{L^2}+\|\nabla m_k\|_{H^1}^2)\leq \|m\|^2_{L^2}+\|\nabla m\|_{H^1}^2=\|m\|^2_{H^1}.$$
    By \eqref{eq002} it follows that $m_k\rightarrow m$ strongly in $H^1$, concluding the proof of the claim.
    
\medskip

By Claim 1 we deduce that the functional $v\mapsto d^2I_{f,h}(u_{f,h})(v,v)$ attains a minimum at $m$. It follows that
$$d^2I_{f,h}(u_{f,h})(m,w)=0,\quad\mbox{ for all }w\in H^1(M),$$
which implies that $m\in H^1(M)$ is a weak solution of the equation
\begin{equation}\label{eq003}
    \begin{cases}
    -\Delta_gm=2fe^{2u_{f,h}}m & \mbox{ in }M,\\
    \dfrac{\partial m}{\partial \nu_g}=he^{u_{f,h}}m & \mbox{ on }\partial M.
\end{cases}
\end{equation}

\noindent{\bf Claim 2:} For $m\in H^1(M)$ given by Claim 1, it holds
$$d^4I_{f,h}(u_{f,h})(m,m,m,m)=-8\int_Mfm^4e^{2u_{f,h}}dv_g-\int_{\partial M}hm^4e^{u_{f,h}}d\sigma_g<0.$$

First, note that $m$ is not constant. Otherwise, this constant cannot be zero, since $\|m\|_{H^1}=1$. In this case, \eqref{eq003} implies that $f$ and $h$ are identically zero. But the Gauss-Bonnet Theorem implies
    $$\int_Mfe^{2u_{f,h}}dv_g+\int_{\partial M}he^{u_{f,h}}d\sigma_g=\mathcal X(M)<0,$$
    which is a contraction. By the first equation in \eqref{eq003} we obtain
    $$2fe^{2u_{f,h}}m^4=-m^3\Delta_gm=-\frac{1}{4}\Delta_gm^4+3|\nabla m|^2_gm^2.$$
    Thus, using that $h\leq 0$, the second equation in \eqref{eq003} and integrating by parts, we get
    \begin{align*}
        d^4I_{f,h}(u_{f,h})&(m,m,m,m)=-8\int_Mfm^4e^{2u_{f,h}}dv_g-\int_{\partial M}hm^4e^{u_{f,h}}d\sigma_g\\
        & = 4\int_Mm^3\Delta_gmdv_g-\int_{\partial M}hm^4e^{u_{f,h}}d\sigma_g\\
        & = \int_M(\Delta_gm^4-12|\nabla m|_g^2m^2)dv_g-\int_{\partial M}hm^4e^{u_{f,h}}d\sigma_g\\
        & = -12 \int_M|\nabla m|_g^2m^2dv_g+\int_{\partial M}m^3\left(4\frac{\partial m}{\partial\nu_g}- hme^{u_{f,h}}\right)d\sigma_g\\
        & = -12 \int_M|\nabla m|_g^2m^2dv_g+\int_{\partial M}3he^{u_{f,h}}m^4\sigma_g<0.
    \end{align*}

Since $u_{f,h}$ is a local minimizer and $d^2I_{f,h}(u_{f,h})(m,m)=0$, then $d^3I_{f,h}(u_{f,h})(m,m,m)=0$. Also we obtain that  $dI_{f,h}(u_{f,h})=0.$ From this, we obtain the expansion
    $$I_{f,h}(u_{f,h}+t m)=I_{f,h}(u_{f,h})+\frac{t^4}{4}d^4I_{f,h}(u_{f,h})(m,m,m,m)+O(t^5)<I_{f,h}(u_{f,h}),$$
    for small $t>0$, where we used  Claim 2. But this leads to a contradiction, since $u_{f,h}$ is a local minimizer.
\end{proof}

\subsection{Existence of strict local minimizers}
We conclude this section with the result that completes the proof of Theorem \ref{Theo_exist}. Following the idea in \cite{MR3628322} (see also \cite{MR1753266}) we prove the next proposition.
\begin{proposition}
Let $(M,g)$ be a compact Riemannian surface with $\partial M\not=\emptyset$ and $\chi(M)<0$. Suppose that for some $f\in C^{0,\alpha}(M)$ and $h\in C^{1,\alpha}(\partial M)$, with $h\leq 0$, the functional $I_{f,h}$ admits a local minimizer $u_{f,h}\in H^1(M)$. Then there exist a neighborhood $U\times W\subset C^{0,\alpha}(M)\times C^{1,\alpha}(\partial M)$ of $(f,h)$  and a smooth map $\mathcal T:U\times W\to H^1(M)$ such that for any $(z,k)\in U\times W$ the function $\mathcal T(z,k)$ is a strict local minimizer of $I_{z,k}$.
\end{proposition}
\begin{proof}

First, we prove the following claim.

\medskip

\noindent{\bf Claim:} There exist two neighborhoods, $U\times W\subset C^{0,\alpha}(M)\times C^{1,\alpha}(\partial M)$ of $(f,h)$  and $W'$ of $(0,0)\in C^{0,\alpha}(M)\times C^{1,\alpha}(\partial M)$, such that for any $(m,z)\in U\times W$ there exists a function $u_{m,z},$ which depends smoothly on $(m,z)$ and solves the boundary value problem:
\begin{equation}\label{eq005}
    \begin{cases}
    -\Delta_gu_{m,z}-1=me^{2u_{m,z}} & \mbox{ in }M,\\
\dfrac{\partial u_{m,z}}{\partial\nu_g}=ze^{u_{m,z}} & \mbox{ on }\partial M.
\end{cases}
\end{equation}

\medskip

Consider the  map $\mathcal Z:C^{2,\alpha}(M)\times C^{0,\alpha}(M)\times C^{1,\alpha}(\partial M)\rightarrow C^{0,\alpha}(M)\times C^{1,\alpha}(\partial M)$
    given by
    $$\mathcal Z(w,m,z)=\left(-\Delta_gw+(f-me^{2w})e^{2u_{f,h}},\frac{\partial w}{\partial\nu_g}+(h-ze^{w})e^{u_{f,h}}\right).$$
    Note that $\mathcal Z(0,f,h)=(0,0),$ the map  $\mathcal Z$ is $C^1$  and its partial derivative with respect to $w$, $D_w\mathcal Z(0,f,h):C^{2,\alpha}(M)\rightarrow C^{0,\alpha}(M)\times C^{1,\alpha}(\partial M)$, is given by
    $$D_w\mathcal Z(0,f,h)\cdot(k)=\left(-\Delta_gk-2fke^{2u_{f,h}},\frac{\partial k}{\partial\nu_g}-hke^{u_{f,h}}\right).$$

    By Proposition \ref{propo001}, Lax-Milgram Theorem and standard elliptic regularity results, we obtain that $D_w\mathcal Z(0,f,h)$ is an isomorphism. By Implicit Function Theorem there is a neighborhood $U\times W\subset C^{0,\alpha}(M)\times C^{1,\alpha}(\partial M)$ of $(f,h)$, a neighborhood $W'\subset  C^{0,\alpha}(M)\times C^{1,\alpha}(\partial M)$ of $(0,0)$, and a smooth map $\mathcal Z_0:U\times W\to W'$ such that for any $(m,z)\in U\times W$ it holds
    $$\mathcal Z(\mathcal Z_0(m,z),m,z)=0.$$
    Since $u_{f,h}$ solves \eqref{eq001}, setting $u_{m,z}=\mathcal Z_0(m,z)+u_{f,h},$ we arrive at the desired  claim.

\medskip

Since $u_{m,z}$ solves \eqref{eq005}, it is a critical point of the functional $I_{m,z}$. The second variation is given by \eqref{eq047}.
Additionally, we assume that $u_{f,h}$ is a local minimizer of the functional $I_{f,h}$, then by Proposition \ref{propo001} we have $d^2I_{f,h}(u_{f,h})(w,w)\geq c_0>0$ for all $w\in H^1(M)$ with $\|w\|^2_{H^1}=1$. Since $u_{m,z}$ depends smoothly on $(m,z)$, we have the convergence $u_{m,z}\rightarrow u_{f,h}$ in $C^{2,\alpha}(M)$ when $(m,z)\rightarrow (f,h)$ in $C^{0,\alpha}(M)\times C^{1,\alpha}(\partial M)$. Choosing $U\times W$  to be smaller if necessary, we have
$$d^2I_{m,z}(u_{m,z})(w,w)\geq c_1>0,$$ 
    for all $(m,z)\in U\times W$ and $w\in H^1(M)$ with $\|w\|^2_{H^1}=1$. Given that $u_{m,z}$ is a critical point of $I_{m,z}$, the Taylor expansion yields that $u_{m,z}$ is a strict local minimizer for $I_{m,z}$.
\end{proof}

\section{Existence of a further critical point of saddle-type}\label{sec:existence-of-a-further}

Let $f\in C^\infty(M)$ and $h\in C^\infty(\partial M)$ be nonconstant functions such that $\max f=\max h=0$. Consider the perturbations $f_\mu=f+\mu$ and $h_{\lambda} = h+\lambda,$ where $\mu$ and $\lambda$ are real numbers. Set $I_{\mu,\lambda}(u):=I_{f_\mu,h_\lambda}(u), u \in H^1\left(M\right)$. As a consequence of Theorem \ref{Theo_exist}, there are $\mu_0>0$ and $\lambda_0>0$ such that for any $\mu\in(0,\mu_0)$ and $\lambda\in(0,\lambda_0)$, the functional $I_{\mu,\lambda}$ admits a smooth strict local minimizer, denoted by $u_{\mu,\lambda}$, which depends smoothly on $\mu$ and $\lambda$. Thus,  as $\mu \downarrow 0$ and $\lambda \downarrow 0$, $u_{\mu,\lambda} \rightarrow u_0$ smoothly in $H^1\left(M \right)$, where $u_0$ is the unique  solution of \eqref{eq001}.

Since $\displaystyle\max f=0$, then given $\mu>0$ there exists a point $p_\mu\in M\backslash\partial M$ such that $f_\mu(p_\mu)>0$. This implies that $I_{\mu,\lambda}$ is no longer bounded from below, independent of $\lambda$. Indeed, let $U\subset M\backslash\partial M$ be a neighborhood of $p_\mu$ such that $f_\mu(p)>0$ for all $p\in U$. Let $w$ be a smooth nonnegative function supported in $U$ and with $w(p_\mu)>0$. We can see that $I_{\mu,\lambda}(tw) \to -\infty$ as $t \to \infty$. Therefore, by Theorem \ref{Theo_exist}, $I_{\mu,\lambda}$ has a mountain pass geometry for all $\mu$ and $\lambda$ small enough. Thus we may expect the existence of another critical point of saddle-type. 

Given that $u_{\mu,\lambda}\to u_0$ as  $(\mu,\lambda)\to (0,0)$, after replacing $\mu_0>0$ and $\lambda_0>0$ with smaller numbers $\mu_0>0$ and $\lambda_0>0$, if necessary, we can find $\rho>0$ such that
\begin{equation}\label{eq014}
I_{\mu,\lambda}\left(u_{\mu,\lambda}\right)=\inf _{\left\|u-u_0\right\|_{H^1}<\rho} I_{\mu,\lambda}(u)  \leq \sup _{\gamma, \nu \in \left(0, \mu_0\right];\theta, \xi \in \left(0, \lambda_0\right]} I_{\gamma,\theta}\left(u_{\nu,\xi}\right)<\beta_{u_0},
\end{equation}
uniformly for all $\mu \in \left(0, \mu_0\right]$ and $\lambda \in \left(0, \lambda_0\right]$, where  
\begin{equation}\label{beta}
    \beta_{u_0}=\inf_{(\mu,\lambda) \in (0,\mu_0]\times(0,\lambda_0] }\{ I_{\mu,\lambda}(u): \rho / 2<\left\|u-u_0\right\|_{H^1}<\rho\}.
\end{equation}
From now on, $\mu_0>0$ and $\lambda_0>0$ are fixed such that $\eqref{eq014}$ and \eqref{beta} are true. 

Given that the functional $I_{\mu,\lambda}$ is unbounded from below, for each $\mu \in  (0,\mu_0]$ and $\lambda \in  (0,\lambda_0]$, consider  $v_{\mu,\lambda} \in H^1\left(M\right)$ such that
\begin{equation}\label{eq015}
    I_{\mu,\lambda}\left(v_{\mu,\lambda}\right)<I_{\mu,\lambda}\left(u_{\mu,\lambda}\right).
\end{equation}
In this case, setting
\begin{equation}\label{eq020}
\Gamma=\left\{p \in C^0\left([0,1] ; H^1\left(M \right)\right): p(0)=u_0\mbox{ and } p(1)=v_{\mu,\lambda}\right\}\footnote{Since $u_{\mu,\lambda} \rightarrow u_0$ as $\mu \downarrow 0$ and $\lambda \downarrow 0$,  we can fix the initial point of comparison paths $p \in \Gamma$  to be $u_0$ instead of $u_{\mu,\lambda},$ provided $\mu_0$ and $\lambda_0$ are sufficiently small.
} ,    
\end{equation}
we have 
\begin{equation}\label{eq007}
    c_{\mu,\lambda}=\inf _{p \in \Gamma} \max _{t \in[0,1]} I_{\mu,\lambda}(p(t)) \geq \beta_{u_0}>I_{\mu,\lambda}\left(u_{\mu,\lambda}\right).
\end{equation}

Note that for every $\mu,\lambda,\gamma\in\mathbb R$ and every $u\in H^1(M)$ it holds
\begin{equation}\label{eq008}
    I_{\mu,\lambda}(u)-I_{\mu,\gamma}(u)=-(\lambda-\gamma)\int_{\partial M}e^ud\sigma_{g}
\end{equation}
and
\begin{equation}\label{eq009}
    I_{\mu,\lambda}(u)-I_{\gamma,\lambda}(u)=-\frac{\mu-\gamma}{2}\int_{ M}e^{2u}dv_{g}.
\end{equation}
Then it follows from Lemma \ref{lemma_bound} that the functions $\mu\mapsto c_{\mu,\lambda}$ and $\lambda\mapsto c_{\mu,\lambda}$ are non-increasing, and therefore differentiable at almost every $\mu\in(0,\mu_0)$ and $\lambda\in(0,\lambda_0)$, respectively. Also, the map $t\mapsto c_{\mu+t,\lambda+t}$ is non-increasing and differentiable at almost every $t$ in its domain. Therefore, the derivative at $t=0$ of the function 
\begin{equation}\label{eq024}
    t\mapsto c_{\mu,\lambda}(t):=c_{\mu+t,\lambda+t}
\end{equation}
exists for almost every $(\mu,\lambda)\in (0,\mu_0)\times(0,\lambda_0)$.

\begin{proposition}\label{propo002}
    Assume the derivative $c_{\mu,\lambda}'(0)$ exists for some $(\mu,\lambda)\in(0,\mu_0)\times(0,\lambda_0)$. Then there exist two sequences, $\{p_n\}$ in $\Gamma$ and $\{t_n\}$ 
in $[0,1]$, such that
\begin{itemize}
    \item[(a)] $I_{\mu,\lambda}\left(p_n(t_n)\right) \rightarrow c_{\mu,\lambda},$
    \item[(b)] $ \displaystyle\max_{0 \leq t \leq 1} I_{\mu,\lambda}\left(p_n(t)\right)\rightarrow c_{\mu,\lambda},$
    \item[(c)] $\left\|d I_{\mu,\lambda}\left(p_n(t_n)\right)\right\| \rightarrow 0$ as $n \rightarrow \infty$.
\end{itemize}
Moreover, $\{p_n(t_n)\}$ satisfies
\begin{equation}\label{eq018}
 \frac{1}{2}\int_{ M} e^{2p_n(t_n)} d v_{g}+\int_{\partial M} e^{p_n(t_n)} d \sigma_{g}\leq -c_{\mu,\lambda}'(0)+3.
\end{equation}
\end{proposition}
\begin{proof}
Consider the decreasing sequences $\mu_n=\mu+1/n$ and $\lambda_n=\lambda+1/n$. Since $c_{\mu,\lambda}$ is min-max, see \eqref{eq007}, we can find two sequences $p_n \in \Gamma$ and $t_n\in[0,1]$ such that
\begin{equation}\label{0001}
\max _{t \in[0,1]} I_{\mu,\lambda}\left(p_n(t)\right) \leq c_{\mu,\lambda}+\frac{1}{n}
\end{equation}
and
\begin{equation}\label{0002}
c_{\mu_n,\lambda_n}-\frac{1}{n} < I_{\mu_n,\lambda_n}(u_n),
\end{equation}
where $u_n=p_n\left(t_n\right)$. By \eqref{eq008}, \eqref{eq009}, \eqref{0001} and \eqref{0002}, we get
\begin{align*}
  0  \leq \frac{1}{2}\int_{M}e^{2u_n}dv_{g}+\int_{\partial M}e^{u_n}d\sigma_{g} & =\frac{I_{\mu,\lambda}(u_n)-I_{\mu_n,\lambda_n}(u_n)}{1/n}  \leq \frac{c_{\mu,\lambda}(0)-c_{\mu,\lambda}(1/n)}{1/n}+2,
\end{align*}
which implies that
\begin{equation}\label{eq011}
    \frac{1}{2}\int_{M}e^{2u_n}dv_{g}+\int_{\partial M}e^{u_n}d\sigma_{g}\leq -c_{\mu,\lambda}'(0)+3
\end{equation}
for $n\in\mathbb N$ sufficiently large. Since  $c_{\mu,\lambda}(t)$ is non-increasing, it follows that  $c_{\mu,\lambda}'(0)\leq 0$. By \eqref{eq011} and Jensen's inequality, it holds
\begin{equation}\label{eq010}
    \fint_M u_ndv_g\leq \log \left(\fint_Me^{2u_n}dv_g\right)\leq C(\mu,\lambda)
\end{equation}
and
$$\fint_{\partial M} u_nd\sigma_g\leq \log\left(\fint_{\partial M}e^{u_n}d\sigma_g\right)\leq C(\mu,\lambda),$$
for all $n\in\mathbb N$ sufficiently large and for some  constant $C(\mu,\lambda)>0$. Here, $\fint$ denotes the mean value. By \eqref{0001} and \eqref{0002}, we obtain
\begin{equation}\label{eq059}
    c_{\mu_n,\lambda_n}-\frac{1}{n}\leq I_{\mu_n,\lambda_n}(u_n)\leq I_{\mu,\lambda}(u_n)\leq \max_{t\in[0,1]}I_{\mu,\lambda}(p_n(t))\leq c_{\mu,\lambda}+\frac{1}{n}.
\end{equation}
This proves items (a) and (b).

We use  \eqref{functional}, \eqref{0001}, \eqref{eq011} and \eqref{eq010} to conclude:
\begin{align}
    \|\nabla u_n\|_{L^2}^2 & =2I_{\mu,\lambda}(u_n)+2\int_Mu_ndv_g+\int_M(f+\mu)e^{2u_n}dv_g+\int_{\partial M}(h+\lambda)e^{u_n}d\sigma_g\label{eq061}\\
    & \leq  2c_{\mu,\lambda}+\frac{2}{n}+C_1\leq C_2,\nonumber
\end{align}
for some positive constant $C_2=C_2(\mu,\lambda)$ independent of $n$. In addition, writing \eqref{eq061}, we obtain
\begin{align*}
    \|\nabla u_n\|_{L^2}^2 -2\int_Mu_ndv_g=2I_{\mu,\lambda}(u_n)+\int_M(f+\mu)e^{2u_n}dv_g+\int_{\partial M}(h+\lambda)e^{u_n}d\sigma_g\leq C_2.
\end{align*}
Together with \eqref{eq010}, we see that the average of $u_n$ is bounded. Therefore, it holds
\begin{equation}\label{eq012}
    \|u_n\|_{H^1}\leq C_3,
\end{equation}
for some positive constant $C_3=C_3(\mu,\lambda,g)$ independent of $n$.

\medskip

First, to prove item (c), prove the following two claims.

\medskip

\noindent{\bf Claim 1:} There exists a constant $C>0$, such that for all $u\in H^1(M)$ satisfying $\|u\|_{H^1}\leq C_3$, the following inequality holds:
$$\|dI_{\mu_n,\lambda_n}(u)-dI_{\mu,\lambda}(u)\|_{H^{-1}}\leq C(|\mu_n-\mu|+|\lambda_n-\lambda|),$$
for all $n\in\mathbb N$.

\medskip

    For every $v\in H^1(M)$, by Hölder's inequality we have
    \begin{align*}
        \langle dI_{{\mu_n,\lambda_n}}(u) & -dI_{\mu,\lambda}(u),v\rangle  =-(\mu_n-\mu)\int_Me^{2u}vdv_g-(\lambda_n-\lambda)\int_{\partial M}e^uvd\sigma_{g} \\
        & \leq |\mu_n-\mu|\left(\int_Me^{4u}dv_g\right)^{1/2}\|v\|_{L^2(M)} + |\lambda_n-\lambda|\left(\int_{\partial M}e^{2u}d\sigma_{g}\right)^{1/2}\|v\|_{L^2(\partial M)}.
    \end{align*}
By a Moser-Trundinger type inequality, see for example \cite{MR925123} and \cite{MR2178789}, and the fact that the trace operator $u\mapsto u|_{\partial M}$ is continuous from $H^1(M)$ to $L^2(\partial M)$, see \cite[Theorem 6.2]{MR3014461}, 
we obtain Claim 1.

\medskip

\noindent{\bf Claim 2:} Given $u,w\in H^1(M)$ with $\|w\|_{H^1}\leq 1$, we have
     $$
\left|I_{\mu,\lambda}(u+w)- I_{\mu,\lambda}(u)-\langle dI_{\mu,\lambda}(u),w \rangle_{H^{-1}\times H^1}\right|\leq C\|w\|^2_{H^1},$$
for some positive constant $C$ independent of $\mu\in(0,\mu_0)$ and $\lambda\in(0,\lambda_0)$.
     
\medskip

By Taylor's expansion, given $x\in M$ there exists $t(x),s(x)\in(0,1)$ such that
\begin{align*}
    & I_{\mu,\lambda}(u+w)  - I_{\mu,\lambda}(u)-\langle dI_{\mu,\lambda}(u),w \rangle_{H^{-1}\times H^1}\\
    & = \frac{1}{2}\int_M|\nabla w|^2dv_g-\frac{1}{2}\int_{M}f_{\mu}(e^{2(u+w)}-e^{2u}-2e^{2u}w)dv_g -\int_{\partial M}h_{\lambda}(e^{u+w}-e^{u}-e^{u}w)d\sigma_g\\
    & \leq \frac{1}{2}\|w\|^2_{H^1}-\int_Mf_{\mu}  e^{2(u+t w)}w^2dv_g-\frac{1}{2}\int_{\partial M}h_{\lambda}e^{u+sw}w^2d\sigma_g\\
    & \leq \frac{1}{2}\|w\|^2_{H^1}+\|f_{\mu}\|_{L^\infty}\int_Me^{2(u+t w)}w^2dv_g+\|h_{\lambda}\|_{L^\infty}\int_{\partial M}e^{u+sw}w^2d\sigma_g.
\end{align*}
By Hölder's inequality and Sobolev's embedding,  we obtain
\begin{align*}
    \int_Me^{2(u+t w)}w^2dv_g \leq\left(\int_Me^{4(u+tw)}dv_g\right)^{1/2}\|w\|^2_{L^4}\leq \left(\int_Me^{8u}dv_g\int_Me^{8|w|}dv_g\right)^{1/4}\|w\|^2_{H^1}.
\end{align*}
Since the trace operator $v\mapsto v|_{\partial M}$ is continuous from $H^1(M)$ to $L^4(\partial M)$, see \cite[Theorem 6.2]{MR3014461}, we get
\begin{align*}
    \int_{\partial M}e^{u+sw}w^2d\sigma_g & \leq\left(\int_{\partial M}e^{2(u+tw)}d\sigma_g\right)^{1/2}\left(\int_{\partial M}w^4d\sigma_g\right)^{1/2}\\
    & \leq C\left(\int_{\partial M}e^{4u}d\sigma_g\int_{\partial M}e^{4|w|}d\sigma_g\right)^{1/4}\|w\|^2_{H^1}.
\end{align*}
Now,  Claim 2 follows from the Moser-Trudinger inequality, see \cite{MR2178789}.

\medskip

Now, assume that there exists $\delta>0$ such that $\|dI_{\mu,\lambda}(p_n(t))\|\geq 2\delta$ for all $n\in\mathbb N$ and all $t\in[0,1]$ with $c_{\mu_n,\lambda_n}-1/n\leq I_{\mu_n,\lambda_n}(p_n(t))$. 

Choose a smooth function $\phi:\mathbb R\rightarrow\mathbb R$ such that $0\leq\phi\leq 1$, $\phi(s)=1$ for $s\geq-1/2$, $\phi(s)=0$ for $s\leq-1$ and $0<\phi(s)<1$ for all $s\in(-1,-1/2)$. For any $n\in\mathbb N$ and  $w\in H^1(M)$ define
$$\phi_n(w):=\phi\left(\frac{I_{\mu_n,\lambda_n}(w)-c_{\mu_n,\lambda_n}}{1/n}\right)$$
and
\begin{equation}\label{eq048}
    \tilde p_n(t):=p_n(t)-\frac{1}{\sqrt{n}}\phi_n(p_n(t))\frac{dI_{\mu,\lambda}(p_n(t))}{\|dI_{\mu,\lambda}(p_n(t))\|},\quad 0\leq t\leq 1.
\end{equation}
By \eqref{eq014}, \eqref{eq015} and \eqref{eq007}, we obtain $$I_{\mu_n,\lambda_n}(p_n(1))\leq I_{\mu_n,\lambda_n}(p_n(0))\leq c_{\mu_n,\lambda_n}-1/n,$$
for $n$ large enough. This implies $\phi_n(p_n(1))=\phi_n(p_n(0))=0$. Thus, $\tilde p_n(1)=p_n(1)$ and $\tilde p_n(0)=p_n(0)$, which implies that $\tilde p_n\in \Gamma$. Fix $s_n\in[0,1]$ with $I_{\mu_n,\lambda_n}(p_n(s_n))\geq c_{\mu_n,\lambda_n}-1/n$ and define $u_n=p_n(s_n)$ and $\tilde u_n=\tilde p_n(s_n)=u_n-\frac{1}{\sqrt{n}}\phi_n(p_n(t))\frac{dI_{\mu,\lambda}(p_n(t))}{\|dI_{\mu,\lambda}(p_n(t))\|}$. Notice that for such $s_n$, we have the inequality \eqref{0002}. In addition, we have inequalities similar to \eqref{eq059} and \eqref{eq012} for $p_n(s_n)$.
 Using the inequality $ab\leq a^2/2+b^2/2$ we obtain
\begin{align*}
    \langle dI_{\mu_n,\lambda_n}(u_n), dI_{\mu,\lambda}(u_n)\rangle 
    & =\|dI_{\mu,\lambda}(u_n)\|^2-\langle dI_{\mu,\lambda}(u_n)-dI_{\mu_n,\lambda_n}(u_n), dI_{\mu,\lambda}(u_n)\rangle\nonumber \\
    & \geq \|dI_{\mu,\lambda}(u_n)\|^2- \|dI_{\mu,\lambda}(u_n)-dI_{\mu_n,\lambda_n}(u_n)\|\| dI_{\mu,\lambda}(u_n)\|\nonumber\\
    &\geq \frac{1}{2}\|dI_{\mu,\lambda}(u_n)\|^2-\frac{1}{2}\|dI_{\mu,\lambda}(u_n)-dI_{\mu_n,\lambda_n}(u_n)\|^2.\nonumber
\end{align*}
By Claim 1 we get
\begin{equation}\label{eq013}
    \langle dI_{\mu_n,\lambda_n}(u_n), dI_{\mu,\lambda}(u_n)\rangle \geq \frac{1}{2}\|dI_{\mu,\lambda}(u_n)\|^2-\frac{1}{2n^2}
\end{equation}
for $n$ large enough. Since $0\leq\phi_n\leq 1$, then $\|p_n(t)-\tilde p_n(t)\|_{H^1}\leq 1$, for all $t\in[0,1]$. By Claim 2, \eqref{eq013} and the inequality $\|dI_{\mu,\lambda}(u_n)\|\geq 2\delta$, we have
\begin{align}
    I_{\mu_n,\lambda_n}(\tilde u_n) & \leq I_{\mu_n,\lambda_n}(u_n)-\frac{\phi_n(u_n)}{\sqrt{n}\|dI_{\mu,\lambda}(u_n)\|}\langle dI_{\mu_n,\lambda_n}(u_n),dI_{\mu,\lambda}(u_n)\rangle+C\frac{\phi_n(u_n)^2}{n}\nonumber\\
    & \leq I_{\mu_n,\lambda_n}(u_n)-\frac{1}{2}\frac{\phi_n(u_n)}{\sqrt{n}}\|dI_{\mu,\lambda}(u_n)\|+\frac{\phi_n(u_n)}{2n^2\sqrt{n}\|dI_{\mu,\lambda}(u_n)\|}+C\frac{\phi_n(u_n)^2}{n}\nonumber\\
    & \leq I_{\mu_n,\lambda_n}(u_n)-\delta\frac{\phi_n(u_n)}{\sqrt{n}}+\frac{\phi_n(u_n)}{4n^2\sqrt{n}\delta}+C\frac{\phi_n(u_n)^2}{n}\nonumber\\
&     \leq I_{\mu_n,\lambda_n}(u_n)-\frac{\delta}{2\sqrt{n}}\phi_n(u_n)\leq  I_{\mu_n,\lambda_n}(u_n)\label{eq060},
\end{align}
for $n$ sufficiently large. Since $\tilde p_n(t)=p_n(t)$ for all $t\in[0,1]$ with $I_{\mu_n,\lambda_n}(p_n(t))\leq c_{\mu_n,\lambda_n}-1/n$, we obtain $I_{\mu_n,\lambda_n}(\tilde p_n(t))\leq I_{\mu_n,\lambda_n}(p_n(t))$, for all $t\in[0,1]$.

Now, notice that if $I_{\mu_n,\lambda_n}(\tilde p_n(t))<c_{\mu_n,\lambda_n}-1/2n$ for all $t\in[0,1]$, then this leads to a contradiction based on the definition of  $c_{\mu_n,\lambda_n}$. This implies that the maximum of $I_{\mu_n,\lambda_n}(\tilde p_n(t))$ can only be achieved at points $t\in[0,1]$ such that  $c_{\mu_n,\lambda_n}-1/2n\leq I_{\mu_n,\lambda_n}(\tilde p_n(t))\leq I_{\mu_n,\lambda_n}(p_n(t))$. From this \eqref{eq007}, \eqref{0001}, \eqref{eq060} and the facts that $I_{\mu_n,\lambda_n}(u)\leq I_{\mu,\lambda}(u)$ and $\phi_n(p_n(t))=1$ if $I_{\mu_n,\lambda_n}(p_n(t))\geq c_{\mu_n,\lambda_n}-1/2n$, we get
\begin{align*}
    c_{\mu_n,\lambda_n} & \leq\max_{t\in[0,1]}I_{\mu_n,\lambda_n}(\tilde p_n(t))\leq\max_{t\in[0,1]:I_{\mu_n,\lambda_n}(p_n(t))\geq c_{\mu_n,\lambda_n}-1/2n}I_{\mu_n,\lambda_n}(\tilde p_n(t))\\
    c_{\mu_n,\lambda_n} & \leq\max_{t\in[0,1]:I_{\mu_n,\lambda_n}(p_n(t))\geq c_{\mu_n,\lambda_n}-1/2n}\left(I_{\mu_n,\lambda_n}(p_n(t))-\frac{\delta}{2\sqrt{n}}\phi_n(p_n(t))\right)\\
& \leq\max_{t\in[0,1]}I_{\mu_n,\lambda_n}(p_n(t))-\frac{\delta}{2\sqrt{n}} \\
& \leq c_{\mu,\lambda}+\frac{1}{n}-\frac{\delta}{2\sqrt{n}}.
\end{align*}
Finally, since $c_{\mu,\lambda}'(0)$ exists and is non positive, it follows that 
$$\frac{c_{\mu,\lambda}(0)-c_{\mu,\lambda}(1/n)}{1/n}\leq -c_{\mu,\lambda}'(0)+1,$$
for all $n$ sufficiently large. Thus,
$c_{\mu,\lambda}\leq c_{\mu_n,\lambda_n}+(-c_{\mu,\lambda}'(0)+1)/n$. Therefore,
$$c_{\mu_n,\lambda_n}\leq c_{\mu_n,\lambda_n}+\frac{1}{\sqrt{n}}\left(\frac{-c_{\mu,\lambda}'(0)+2}{\sqrt{n}}-\frac{\delta}{2}\right)<c_{\mu_n,\lambda_n},$$
for $n$ large enough, which is a contradiction. 

To complete the proof of item (c), choose a sequence $t_n\in[0,1]$ satisfying $c_{\mu_n,\lambda_n}-1/n<I_{\mu_n,\lambda_n}(p_n(t_n))$ such that $\|dI_{\mu,\lambda}(p_n(t_n))\|\to 0$ as $n\to\infty$.  This finishes the proof of theorem.
\end{proof}

\subsection{Proof of Theorem \ref{second_critical}}
In this subsection we complete the Proof of Theorem \ref{second_critical}.

\begin{theorem}\label{large_solution}
    Suppose the derivative $c_{\mu,\lambda}'(0)$ exists for some $(\mu,\lambda)\in(0,\mu_0)\times(0,\lambda_0)$. Then the functional $I_{\mu,\lambda}$ admits a critical point $u^{\mu,\lambda}$ with energy $I_{\mu,\lambda}(u^{\mu,\lambda})=c_{\mu,\lambda}$ and bounds on area  and length 
    \begin{equation}\label{eq019}
        \frac{1}{2}\int_{M}e^{2u^{\mu,\lambda}}dv_{g}+\int_{\partial M}e^{u^{\mu,\lambda}}d\sigma_{g}\leq -c_{\mu,\lambda}'(0)+3,
    \end{equation}
     and such that $u^{\mu,\lambda}$ is not a local minimizer of $I_{\mu,\lambda}$.
\end{theorem}
\begin{proof}
    The proof of Proposition \ref{propo002} gives us two sequences $\{p_n\}$ in $\Gamma$ and $\{t_n\}$ in $[0,1]$ such that $u_n=p_n(t_n)\in H^1(M)$ satisfies items (a)--(c) of Proposition \ref{propo002}, \eqref{eq018} and \eqref{eq012}. From \eqref{eq012}, we can assume that $u_n \rightharpoonup u^{\mu,\lambda}$ weakly in $H^1\left(M\right)$ and $u_n\rightarrow u^{\mu,\lambda}$ strongly in $L^2(M)$, for some $u^{\mu,\lambda}\in H^1(M)$. Then, up to a subsequence, we can assume  that $e^{2 u_n} \rightarrow e^{2 u^{\mu,\lambda}}$ in $L^2(M),$ and since the trace embedding $T: H^1(M)\to L^2(\partial M)$ is compact, see \cite[Theorem 6.2]{MR3014461}, we have $e^{u_n}\rightarrow e^{u^{\mu,\lambda}}$ in $L^2(\partial M)$, see, for instance, \cite[Proposition 3.16]{phdthesis} for a proof (See also the proof of \cite[Lemma 4.1]{MR3836128}). From this and \eqref{eq018} we obtain \eqref{eq019}. Note that
\begin{align*}
    \left\|\nabla u_n-\nabla u^{\mu,\lambda}\right\|_{L^2(M)}^2  = \left\langle d I_{\mu,\lambda}\left(u_n\right), u_n-u^{\mu,\lambda}\right\rangle-\int_M\left\langle\nabla u^{\mu,\lambda}, \nabla u_n-\nabla u^{\mu,\lambda}\right\rangle dv_{g}\\
     +\int_M\left(u_n-u^{\mu,\lambda}\right) d v_{g}  +\int_M f_\mu e^{2 u_n}\left(u_n-u^{\mu,\lambda}\right) dv_{g}+\int_{\partial M} h_\lambda e^{u_n}\left(u_n-u^{\mu,\lambda}\right) d \sigma_{g}.
\end{align*}
This implies that $\nabla u_n\rightarrow\nabla u^{\mu,\lambda}$ in $L^2(M)$.
Thus, we have that $u_n \rightarrow u^{\mu,\lambda}$ strongly in $H^1\left(M\right)$. Consequently  $I_{\mu,\lambda}\left(u_n\right) \rightarrow I_{\mu,\lambda}\left(u^{\mu,\lambda}\right)$ and $d I_{\mu,\lambda}\left(u_n\right) \rightarrow d I_{\mu,\lambda}\left(u^{\mu,\lambda}\right)$ as $n \rightarrow \infty$, and by Proposition \ref{propo002} the function $u^{\mu,\lambda}$ is a critical point for $I_{\mu,\lambda}$ at level $I_{\mu,\lambda}\left(u^{\mu,\lambda}\right)=c_{\mu,\lambda}$.

Reasoning by contradiction, assume  that $u^{\mu,\lambda}$ is a local minimizer of $I_{\mu,\lambda}$ for $\mu>0$ and $\lambda>0$ small enough. By Theorem \ref{Theo_exist}, the function $u^{\mu,\lambda}$ is a strict local minimizer. The sequences $\{p_n\}$ and $(u_n)$, as before, satisfy $\displaystyle\max_{0 \leq t \leq 1} I_{\mu,\lambda}\left(p_n(t)\right)\rightarrow c_{\mu,\lambda}$ and $I_{\mu,\lambda}(u_n)\rightarrow c_{\mu,\lambda}$. Since $u_n\rightarrow u^{\mu,\lambda}$ in $H^1(M)$, there exists  a function $v$  arbitrarily close to $u^{\mu,\lambda}$ such that $I_{\mu,\lambda}(v)$ is less than $I_{\mu,\lambda}(u^{\mu,\lambda})$. This is a contradiction.
\end{proof}

From Theorem \ref{large_solution} we obtain Theorem \ref{second_critical}.

\section{Blow-up analysis and Proof of Theorem \ref{theor-blow-up-analysis}}\label{sec:blow-up-analysis}

In this section, we  assume that the functions $f\in C^\infty(M)$ and $h\in C^\infty(\partial M)$ are nonconstant with $\max f=\max h=0$, the maximum point of $f$ and $h$ are the same, non-degenerate, and critical points of $f$. 

We begin with the following lemma, where we show an explicit estimate of the mountain-pass energy level  $c_{\mu,\lambda}$ associated with  the set 
$\Gamma,$ as defined in \eqref{eq020}. First, for convenience we fix some notations, as follows: The euclidean metric will be denoted by $g_{euc}$. Given $p_0\in\mathbb R^2$ and $r>0,$ define
\begin{align*}
    B_r(p):= & \{q\in\mathbb R^2: |p-q|<r\};\\
    B^+_r(p):= & \{(s,t)\in B_r(p)\subset\mathbb R^2:\:  t\geq 0\};\\
    \partial^+B_r^+(p):= & \{(s,t)\in \partial B_r(p)\subset\mathbb R^2:\:  t> 0\};\\
    \Gamma_r(p):= & \{(s,t)\in \partial B^+_r(p):\:  t=0\}.
\end{align*}
These sets can also be denoted by $B_r$, $B_r^+$, $\partial^+B_r^+$ and $\Gamma_r$. 

Before starting, remember that we are assuming $K_g=-1$ and $\kappa_g=0$.

\begin{lemma}\label{lemma_bound}
For any constant $\kappa>2\pi$ there exists $(\mu_\kappa,\lambda_\kappa)\in(0,\mu_0/2]\times(0,\lambda_0/2]$ such that for any $(\mu,\lambda)\in(0,\mu_\kappa)\times(0,\lambda_\kappa)$ there exists $v_{\mu,\lambda}\in H^1(M)$ such that choosing $v_{\gamma,\nu}=v_{\mu,\lambda}$ for every $(\gamma,\nu)\in[\mu,2\mu]\times[\lambda,2\lambda]$ in \eqref{eq020}, there holds 
    $
    I_{ \gamma,\nu}(v_{ \gamma,\nu})=I_{ \gamma,\nu}(v_{ \mu,\lambda})< I_{ \gamma,\nu}(u_{ \gamma,\nu}),
    $ 
    and the number $c_{\gamma,\nu}$ is independent of $\mu$ and $\lambda$ such that $(\gamma,\nu)\in [\mu,2\mu]\times[\lambda,2\lambda]$. 
    Moreover, we obtain the bound $c_{\gamma,\nu}\leq \kappa\log(2/\gamma)$.
\end{lemma}
\begin{proof}

Let $p_0\in \partial M$ be a maximum point of $f$. By assumption, it holds $f(p_0)=h(p_0)=0$. Consider local conformal coordinates $x=(x_1,x_2)$ around $p_0=0$ such that $g=e^{2v_0}g_{euc}$, for some smooth function $v_0$ with $v_0(0)=0$ and $\partial M$ correspond to $x_2=0$. 

Since $p_0$ is a non-degenerate maximum point of $f$ and $h$, we have that 
    $$f(x)=\frac{1}{2}\nabla f(0)(x,x)+O(|x|^3)\geq -\frac{\mu}{2},\quad\mbox{ and } \quad f_\mu(x)=f(x)+\mu\geq\frac{\mu}{2},\quad\forall x\in B_{\sqrt{\mu}/L}^+(0),$$
    and
    $$h(x)=\frac{1}{2}\nabla h(0)(x,x)+O(|x|^3)\geq-\frac{\lambda}{2},\quad\mbox{ and }\quad h_\lambda(x)=h(x)+\lambda\geq\frac{\lambda}{2},\quad\forall x\in \Gamma_{\sqrt{\lambda}/L}(0),$$
     for all $(\mu,\lambda)\in(0,\mu_0]\times(0,\lambda_0]$, and for some constant $L>\max\{\sqrt{\mu_0},\sqrt{\lambda_0}\}$ which does not depend on $\lambda$ and $\mu$. Consider a function $z_\mu\in H^1(B_1^+(0))$ given by 
    $$z_\mu(x)=\begin{cases}
        -\log\mu,&\mbox{if }|x|\leq\mu,\\
        -\log|x|, & \mbox{if }\mu\leq|x|\leq 1.
    \end{cases}$$
In these local conformal coordinates, we define a function $w_\mu\in H^1(M)$ as $w_\mu(x)=z_\mu(Lx/\sqrt{\mu})$ in $B^+_{\sqrt{\mu}/L}(p_0)$ and $w_\mu\equiv 0$ outside of $B^+_{\sqrt{\mu}/L}(p_0)$.    This implies that
    $$\int_M|\nabla_g w_\mu|^2_g  dv_g  =\int_{B^+_{\frac{\sqrt{\mu}}{L}}(p_0)}|\nabla_g w_\mu|^2_gdv_g=\int_{B^+_{\frac{\sqrt{\mu}}{L}}(0)}|\nabla w_\mu(x)|^2dx$$
    $$= \frac{L^2}{\mu}\int_{B^+_{\frac{\sqrt{\mu}}{L}}(0)}|\nabla z_\mu(Lx/\sqrt{\mu})|^2dx =\int_{B^+_1(0)}|\nabla z_\mu(y)|^2dy=\pi\int_\mu^1\frac{1}{r}dr=-\pi\log\mu.$$
Using Young's inequality $2ab\leq \delta a^2+b^2/\delta$, with $\delta=\frac{\kappa-2\pi}{4\pi}$, $a=s\|\nabla_g w_\mu\|^2_{L^2(M)}$ and $b=\|\nabla_g u_0\|^2_{L^2(M)}$, we get
\begin{align*}
    \|\nabla_g(u_0+sw_\mu)\|^2_{L^2(M)} & \leq (1+\delta)s^2\|\nabla_gw_\mu\|^2_{L^2(M)}+(1+1/\delta)\|\nabla_g u_0\|^2_{L^2(M)}\\
    & =-(\kappa+2\pi)\frac{s^2}{4}\log\mu+C,
\end{align*}
where $C=C(u_0,\kappa)>0$. Also, we have
\begin{align*}
    \int_Mf_\mu & e^{2(u_0+sw_\mu)}  dv_g  \geq\frac{\mu}{2}\int_{B_{\frac{\sqrt{\mu}}{L}}^+(p)}e^{2(u_0+sw_\mu)}dv_g-C \geq c\frac{\mu}{2}\int_{B_{\frac{\sqrt{\mu}}{L}}^+(p)}e^{2sw_\mu}dv_g-C\\
    & \geq c\frac{\mu^2}{2L^2}\int_{B_1^+(0)}e^{2sz_\mu(x)}dx-C \geq c\frac{\mu^2}{2L^2}\int_{B_\mu^+(0)}e^{2sz_\mu(x)}dx-C \geq \frac{\pi}{2}cL^{-2}\mu^{4-2s}-C
\end{align*}
and
\begin{align*}
    \int_{\partial M}h_\lambda e^{u_0+sw_\mu}d\sigma_g  \geq \frac{\lambda}{2}\int_{-\sqrt{\mu}/L}^{\sqrt{\mu}/L}e^{u_0+sw_\mu}dt-C\geq -C.
\end{align*}
In the last integral, we identify the functions defined in an interval of $\mathbb R$ with the function defined in a subset of $\partial M$. Hence, using \eqref{functional} and that $w_\mu\geq 0$, we obtain
$$I_{\mu,\lambda}(u_0+sw_\mu)\leq -\frac{\kappa+2\pi}{2}\frac{s^2}{4}\log\mu-\frac{\pi}{2}cL^{-2}\mu^{4-2s}+C.$$
Note that for $\mu,\lambda\in(0,1)$ it holds $I_{\mu,\lambda}(u_0+sw_\mu)\rightarrow-\infty$ as $s\rightarrow+\infty$. By \eqref{eq008} and \eqref{eq009}, we may fix a sufficiently large  $s_{\mu,\lambda}>0$ and define $v_{\mu,\lambda}:=u_0+s_{\mu,\lambda}w_\mu$ (see \eqref{eq015}), satisfying 
$$   
I_{\mu,\lambda}(v_{\mu,\lambda})<\inf_{(\gamma,\nu)\in(0,\mu_0)\times(0,\lambda_0)}I_{\gamma,\nu}(u_{\gamma,\nu}).
$$
Thus, if $p(t):=u_0+ts_{\mu,\lambda}w_\mu$, by \eqref{eq007}, we obtain
\begin{align*}
    c_{\mu,\lambda} & \leq\max_{t\in[0,1]}I_{\mu,\lambda}(p(t)) \leq \sup_{s>0}\left(\frac{\kappa+2\pi}{2}\frac{s^2}{4}\log\frac{1}{\mu}-\frac{\pi}{2}cL^{-2}\mu^{4-2s}+C\right).
\end{align*}

Consider the function $\varphi(s)=\frac{\kappa+2\pi}{2}\frac{s^2}{4}\log\frac{1}{\mu}-\frac{\pi}{2}cL^{-2}\mu^{4-2s}+C$ with derivative $\varphi'(s)=\frac{\kappa+2\pi}{4}\log\frac{1}{\mu}\left(s-\frac{4\pi c}{L^2(\kappa+2\pi)}\mu^{4-2s}\right)$. Since $\varphi(s)\to-\infty$ as $s\to\infty$, then $\varphi$ has a maximum point at some $s(\mu)\in[0,\infty)$. Note that $\varphi'(0)<0$ and $\varphi(2)>\varphi(0)$, for sufficiently small $\mu>0$. This implies that the maximum point of $\varphi$ belongs to $(0,\infty)$. The critical points of $\varphi$ are the solutions of equation $s=\frac{4\pi c}{L^2(\kappa+2\pi)}\mu^{4-2s}$, one of which is the maximum point. This equation has only two solutions for $\mu\in(0,1)$. Since $\varphi'(0)<0$, $\varphi(2)>\varphi(0)$ and $\varphi'(2)>0$, for $L\gg 1$, we have $s(\mu)>2$, otherwise $\varphi$ would have more than two critical points. In addition, the condition $\varphi'(s(\mu))=0$ implies that $s(\mu)\to 2$ as $\mu\to 0^+$.

Thus, for all small enough $\mu>0$ we have
$$ c_ {\mu,\lambda}< \frac{\kappa+2\pi}{2}\log\frac{1}{\mu}\leq \kappa\log\frac{1}{\mu}.$$

Since $I_{\gamma,\nu}(v_{\mu,\lambda})\leq I_{\mu,\lambda}(v_{\mu,\lambda})$ for $\gamma>\mu$  and $\nu>\lambda,$ the same comparison function $v_{\mu,\lambda}$ can be used for every $(\gamma,\nu)\in (\mu,2\mu)\times(\lambda,2\lambda)\subset (0,\mu_0/2]\times(0,\lambda_0/2]$. Recall that, as mentioned above, we can assume that $\mu_0$ and $\lambda_0$ are sufficiently small. For such $(\gamma,\nu)$ we obtain the bound
\begin{equation*}
I_{\gamma,\nu}\left(v_{\mu,\lambda}\right)<I_{\gamma,\nu}\left(u_{\gamma,\nu}\right)<\sup_{(s,t)\in [\mu,2\mu]\times [\lambda,2\lambda]}I_{\gamma,\nu}\left(u_{s,t}\right)<\beta_{u_0}\leq c_{\gamma,\nu}<\kappa \log\frac{1}{\mu}<\kappa\log \frac{2}{\gamma},
\end{equation*}
where $\beta_{u_0}$ is defined in \eqref{beta}. Finally, since $v_{\mu,\lambda}$   depends continuously on $\mu$ and $\lambda$,  by construction, with $I_{\mu,\lambda}(v_{\mu,\lambda}) < \inf\{I_{\gamma,\nu}(u_{\gamma,\nu}):(\gamma,\nu)\in (0,\mu_0/2]\times(0,\lambda_0/2]\},
$
the number $c_{\gamma,\nu}$ is defined independently of $\mu$ and $\lambda$ such that $(\gamma,\nu)\in (\mu,2\mu)\times (\lambda,2\lambda).$
\end{proof}

It is important to mention that Lemma \ref{lemma_bound} holds if $f(p_0)=0$ and $h(p_0)<0$, but since we are not interested in this case, we omitted the proof.

\begin{lemma}\label{lem003}
    For $c_{\mu,\lambda}(t)$ defined as in \eqref{eq024} it holds
    \begin{equation}\label{eq025}
    \liminf_{(\mu,\lambda)\downarrow(0,0)}\mu\left|c_{\mu,\lambda}'(0)\right| \leq 2\pi\quad\mbox{ and }\quad\liminf_{(\mu,\lambda)\downarrow(0,0)}\lambda\left|c_{\mu,\lambda}'(0)\right| \leq 2\pi.
\end{equation}

\end{lemma}
\begin{proof}
        For the first inequality, suppose that for $\kappa>\kappa_1>2\pi,$ we have $|c_{\mu,\lambda}'(0)|\geq \kappa\mu^{-1}$ for almost every $\mu$ and $\lambda$ small enough. Then for any $\lambda_1>\mu_1>\mu_2>0$ sufficiently small, using the fact  that $c_{\mu,\lambda}(t)$ is non-increasing, we obtain
\begin{align*}
    c_{\mu_2,\lambda_1-\mu_1+\mu_2}-c_{\mu_1,\lambda_1} & \geq\int_{\mu_2}^{\mu_1}\left|\frac{d}{dt}c_{t,\lambda_1-\mu_1+t}\right|dt =\int_{\mu_2}^{\mu_1}|c_{t,\lambda_1-\mu_1+t}'(0)|dt\geq \int_{\mu_2}^{\mu_1}\frac{\kappa}{t}dt=\kappa\log\frac{\mu_1}{\mu_2}.
\end{align*}
By  Lemma  \ref{lemma_bound},  we have $c_{\mu_2,\lambda_1-\mu_1+\mu_2}\leq \kappa_1\log(2/ \mu_2)$ for all sufficiently small $\mu_2$. Thus
$$c_{\mu_1,\lambda_1}\leq(\kappa_1-\kappa)\log\frac{2}{\mu_2}-\kappa\log\mu_1+\kappa\log 2,$$
for all small enough $\mu_2>0$, which is a contraction by \eqref{eq007}, since $\kappa_1<\kappa$.

Now, suppose that for $\kappa>\kappa_1>2\pi$ it holds $|c'_{\mu,\lambda}(0)|\geq\kappa\lambda^{-1}$ for almost every $\lambda$ and $\mu$ small enough. Then for any $\mu_1>\lambda_1>\lambda_2>0$ sufficiently small, as before we have
\begin{align*}
    c_{\mu_1-\lambda_1+\lambda_2,\lambda_2}-c_{\mu_1,\lambda_1} & \geq \int_{\lambda_2}^{\lambda_1}\left|\frac{d}{dt}c_{\mu_1-\lambda_1+t,t}\right|dt = \int_{\lambda_2}^{\lambda_1}|c_{\mu_1-\lambda_1+t,t}'(0)|dt\geq\int_{\lambda_2}^{\lambda_1}\frac{\kappa}{t}dt=\kappa\log\frac{\lambda_1}{\lambda_2}.
\end{align*}
Fix $\mu_1$ and choose $\lambda_1=\mu_1-\lambda_2$. For all sufficiently small $\lambda_2>0$, we have  by  Lemma  \ref{lemma_bound}  that $c_{\mu_1-\lambda_1+\lambda_2,\lambda_2}\leq \kappa_1\log (2/(\mu_1-\lambda_1+\lambda_2))\leq-\kappa_1\log(\lambda_2)$.  Thus, 
$$c_{\mu_1,\lambda_1}\leq(\kappa_1-\kappa)\log\frac{1}{\lambda_2}-\kappa\log\lambda_1,$$
for all small enough $\lambda_2>0.$  Again, this leads to a contradiction by \eqref{eq007}, since $\kappa_1<\kappa$.
\end{proof}

\begin{proposition}\label{propo003}
      There exist a constant $C>0$, two sequences $\lambda_n\downarrow 0$ and $\mu_n\downarrow 0$, with $\lambda_n^2-\lambda_n^3\leq\mu_n\leq \lambda_n^2+\lambda_n^3$, and a 
      corresponding solutions $u_n:=u^{\mu_n,\lambda_n}\neq u_{\mu_n,\lambda_n}$ of \eqref{eq017} of “mountain pass” type inducing conformal metrics $g_n = e^{2u_n}g$ with total curvature satisfying
\begin{equation}\label{total_curvature}
  \int_M|K_{g_n}|dv_{g_n}+\int_{\partial M}|\kappa_{g_n}|d\sigma_{g_n} \leq C
  \end{equation}
for some constant $C>0$ independently of $n$, and
\begin{equation}\label{eqbounds}
    \limsup_{n\rightarrow+\infty}\left(\frac{1}{2}\mu_n\int_Me^{2u_n}dv_g+\lambda_n\int_{\partial M}e^{u_n}d\sigma_g\right)\leq 2\pi.
\end{equation}
\end{proposition}
\begin{proof}
Given $\overline{\lambda} > 0$, the function $t \mapsto c_{\overline{\lambda}^2 + t, \overline{\lambda} + t}$ is differentiable in almost every $t$ in its domain. Furthermore, if $\mu = \overline{\lambda}^2 + t$ and $\lambda = \overline{\lambda} + t$, then $\mu \in (\lambda^2 - \lambda^3, \lambda^2 + \lambda^3)$ for all $t \geq 0$ small enough. Thus, we can find two sequences $\lambda_n \downarrow 0$ and $\mu_n \downarrow 0$, with $\mu_n/\lambda_n^2 \in (1 - \lambda_n, 1 + \lambda_n)$ such that $c_{\mu_n, \lambda_n}'(0)$ exists. Thus, for each $n$, Theorem \ref{large_solution} gives us the desired solution $u_n$. By \eqref{eq017} and the Gauss-Bonnet Theorem, it follows that 
$$
\int_M f_{\mu_n}dv_{g_n} +\int_{\partial M} h_{\lambda_n}d\sigma_{g_n} =\int_M K_{g_n}dv_{g_n}+\int_{\partial M} \kappa_{g_n}d\sigma_{g_n} =2 \pi \chi(M).
$$
Therefore
$$
-\int_M f e^{2 u_n} dv_g-\int_{\partial M} h e^{u_n}d\sigma_{g}=\mu_n \int_M e^{2 u_n} dv_{g}+\lambda_n \int_{\partial M} e^{u_n} d \sigma_{g}-2 \pi \chi(M).
$$
 By \eqref{eq019} we have
 $$\frac{1}{2}\int_M e^{2 u_n} dv_{g}+\int_{\partial M} e^{u_n} d \sigma_{g}\leq  |c_{\mu_n,\lambda_n}'(0)|+3. $$ 
 Since $0<\mu_n<\lambda_n^2+\lambda_n^3\leq\lambda_n$, using \eqref{eq025} we obtain \eqref{eqbounds}. Finally, for $n$ large enough such that $\mu_n\leq \lambda_n/2$, it holds
$$
\int_M|f|e^{2u_n}dv_g+\int_{\partial M}|h|e^{u_n}d\sigma_g\leq \lambda_n|c_{\mu_n,\lambda_n}'(0)|+ C\leq C <\infty
$$
uniformly in $n\in\mathbb N$. From this bound, we obtain \eqref{total_curvature}.
\end{proof}

\begin{lemma} \label{lem002}
Let $u_n$ be the sequence of solutions of \eqref{eq017} given by Proposition \ref{propo003}. Consider an open set $\Omega\subset\subset M_-=\{x\in M: f<0\mbox{ and }h<0\}$.
\begin{enumerate}
    \item[(a)] There exists a constant $C>0$ such that for all $n$ it holds
    $u_n>-C$.
    \item[(b)] Let $u_n^+=\max\{u_n,0\}$. Then, for $n$ sufficiently large it holds
$$\int_{\Omega}\left(\left|\nabla u_n^{+}\right|^2+(u_n^+)^2\right)dv_g+\sup_\Omega(-h)\int_{\partial\Omega\cap\partial M}(u_n^+)^2d\sigma_g\leq C(\Omega).$$
\item[(c)] There exists a constant $C(\Omega)$ such that 
$u_n\leq C(\Omega)$, for all sufficiently large $n$.
\end{enumerate}
\end{lemma}
\begin{proof}
First, let us prove item (a). The proof is based on the ideas of \cite{MR1257102} and \cite{xu2023kazdan}. Consider the following problem: 
\begin{equation}\label{eq022}
    \begin{cases}
       -\Delta_gv= b_1, & \mbox{in } M,\\
       \dfrac{\partial v}{\partial \nu_g}=b_2, & \mbox{on }\partial M,
    \end{cases}
\end{equation}
    where $b_1\in(0,1)$ and $b_2=-b_1\frac{\mbox{vol}(M)}{\mbox{vol}(\partial M)}<0$. Note that due to the choices of $b_1$ and $b_2$, \eqref{eq022} has a solution $v$. Thus, there exists a constant $a>0$ such that for all $c\geq a$, the function $\varphi_c=v-c$ satisfies
    $$\begin{cases}
       -\Delta_g\varphi_c-1-f_{\mu_n}e^{2\varphi_c}= b_1-1-f_{\mu_n}e^{2\varphi_c}<0, & \mbox{in } M,\\
       \dfrac{\partial \varphi_c}{\partial \nu_g}-h_{\lambda_n}e^{\varphi_c}=b_2-h_{\lambda_n}e^{\varphi_c}<0, & \mbox{on }\partial M.
    \end{cases}$$

Item (a) will follow if $u_n\geq\varphi_{a}$, for all $n$. 

If there exists some $n$ with $u_n(x)<\varphi_a(x)$, for some $x\in M$, using that $M$ is compact and that $\sup\varphi_c\to-\infty$ when $c\to+\infty$, we can find some $c>a$ such that $u_n\geq\varphi_c$ on $M$, and $u_n(p)=\varphi_c(p)$, for some $p\in M$. By \eqref{eq017} we find that
\begin{equation}\label{eq023}
    \begin{cases}
    -\Delta_g(u_n-\varphi_c)-f_{\mu_n}(e^{2u_n}-e^{2\varphi_c})>0,\\
    \dfrac{\partial}{\partial\nu_g}(u_n-\varphi_c)-h_{\lambda_n}(e^{u_n}-e^{\varphi_c})>0.
\end{cases}
\end{equation}
If $p\not\in\partial M$, by the maximum principle we conclude that $u_n\equiv\varphi_c$, which is impossible. If $p\in\partial M$, Hopf's Lemma implies that
$$\dfrac{\partial}{\partial\nu_g}(u_n-\varphi_c)(p)<0.$$
But, the boundary equation of \eqref{eq023} gives us
$$\dfrac{\partial}{\partial\nu_g}(u_n-\varphi_c)(p)>0,$$
which is a contradiction.

Let us prove item (b). 
Given any domain $\Omega\subset\subset M_-$, it is enough to prove the result for any geodesic ball $B_r^+\subset\subset M_-$, since afterward the estimate for $\Omega$ can be deduced by a covering argument.

 Let $B_{2r} \subset M_{-}$ be the geodesic ball centered at some point  $p\in M$ with radius $4r=\mbox{dist}(p, \partial M_-\backslash\partial M) > 0$. If $B_{2r}\cap \partial M=\emptyset$, then the result follows from the proof of \cite[Lemma 2]{MR1257102}, see also the appendix in \cite{MR3351750}. Thus, suppose that $p\in\partial M$ and denote the geodesic ball as $B_{2r}^+$. 
 
 Let $\psi$ be a smooth cut-off function $0 \leq \psi \leq 1$ supported in $B_{2r}^+$ such that 
 $\psi\equiv 1$ in  $B_{r}^+$ and $|\nabla\psi|\leq Ar^{-1}$, for some positive constant $A$. Multiply the first equation in \eqref{eq017} by $\psi^4u_n^+$ and use integration by parts to obtain
\begin{equation}\label{eq033}
\int_{B_{2r}^+}\left(\left\langle\nabla u_n^{+}, \nabla\left(\psi^4 u_n^{+}\right)\right\rangle-\psi^4 u_n^{+}-f_n e^{2 u_n^{+}} \psi^4 u_n^{+}\right)dv_g-\int_{\Gamma_{2r}}h_{\lambda_n} e^{u_n^{+}} \psi^4 u_n^{+}d\sigma_g=0,    
\end{equation}
where $\Gamma_{2r}=\partial M\cap B_{2r}^+.$ Note that there exists $\varepsilon>0$ such that $f_{\mu_n},h_{\lambda_n} \leq-\varepsilon$, on $B^+_{2r}$ for sufficiently large $n$, and
$\langle \nabla u_n^+,\nabla(\psi^4u_n^+)\rangle=|\nabla(\psi^2u_n^+)|^2-|\nabla\psi^2|^2(u_n^+)^2.$
Using Young's inequality we get
$$|\nabla\psi^2|^2(u_n^+)^2=4|\nabla\psi|^2(\psi u_n^+)^2\leq 4A^2r^{-2}(\psi u_n^+)^2\leq\frac{1}{2}\varepsilon(\psi u_n^+)^4+C,$$
where $C=C(\varepsilon,r,A)$. Also, we estimate
$\psi^4 u_n^{+} \leq \frac{1}{2} \varepsilon \left(\psi u_n^{+}\right)^4+C$. Since $e^{2 t} \geq t^3$ and $e^t\geq t$ for $t \in \mathbb{R}$,  we obtain  from \eqref{eq033} the following inequality
 \begin{align*}
\int_{B_{2r}^+}\left(\left|\nabla\left(\psi^2 u_n^{+}\right)\right|^2 \right.&\left.+\varepsilon \psi^4\left(u_n^{+}\right)^4\right) dv_g+\varepsilon\int_{\Gamma_{2r}}\psi^4(u_n^+)^2\leq\int_{B_{2r}^+}\left(\langle \nabla u_n^+,\nabla(\psi^4u_n^+)\rangle\right.\\
& \left.+|\nabla\psi^2|^2(u_n^+)^2 -f_{\mu_n} e^{2u_n^+}\psi^4u_n^{+}\right)dv_g -\int_{\Gamma_{2r}}h_{\lambda_n} e^{u_n^{+}} \psi^4 u_n^{+}d\sigma_g\\
& =\int_{B_{2r}^+}\left(|\nabla \psi^2|^2\left(u_n^{+}\right)^2+ \psi^4 u_n^{+}\right)dv_g\leq  \varepsilon\int_{B_{2r}^+} \left(\psi u_n^{+}\right)^4dv_g+C,  
 \end{align*}
 which implies
 $$\int_{B_{2r}^+}\left|\nabla\left(\psi^2 u_n^{+}\right)\right|^2dv_g+\varepsilon\int_{\Gamma_{2r}}\psi^4(u_n^+)^2d\sigma_g\leq C,$$
 where $C$ depends on $r$ and $\varepsilon$. Hence the claim follows from Poincaré's inequality.

It remains to prove Item (c). Again, we can assume that $\Omega\cap\partial M\not=\emptyset$. Given $p\in\partial M\cap M_-$, consider a local conformal chart $\Psi:B_{2r}^+(p)\to B_{2r}^+:=B_{2r}^+(0)\subset\mathbb R^2_+$ such that $\overline{B_{2r}^+(p)}\subset M_-$ and $g=e^{2v_0}g_{euc}$, where $v_0$ is a smooth function defined in $B_{2r}^+$. This implies that $v_n:=u_n+v_0$ satisfies the equation
$$\begin{cases}
    -\Delta v_n=f_{\mu_n}e^{2v_n}, & \mbox{ in } B_{2r}^+\\
    \dfrac{\partial v_n}{\partial\nu}=h_{\lambda_n}e^{v_n} & \mbox{ on } \Gamma_{2r}.
\end{cases}$$

Let us prove that $v_n$ is uniformly bounded from above in $B_r^+$.

From the Moser-Trundiger inequalities \cite{MR3836128,phdthesis}, together with the fact that $\left(u_n^{+}\right)$ is uniformly bounded in $H^1(B^+_{2r})$ and in $L^2(B^+_{2r}\cap \partial M)$, it follows that $f_{\mu_n} e^{2 v_n}$ and $h_{\lambda_n} e^{v_n}$ are uniformly bounded in $L^q(B_{2r}^+)$ and in $L^q(\Gamma_{2r})$, respectively, for any $q>1$. Let $w_n$ be the solution of
$$
\begin{cases}
-\Delta w_n=f_{\mu_n} e^{2 v_n}, & \text { in } B^+_{2r},\\
\dfrac{\partial w_n}{\partial \nu}=h_{\lambda_n} e^{v_n}, & \text { on }  \Gamma_{2r},\\
 w_n=0, & \text { on }  \partial^+ B^+_{2r}.
\end{cases}
$$
By standard elliptic estimates imply that 
$$
\|w_n\|_{L^\infty(\overline{B}^+_{r})}\leq C,
$$
for some constant $C>0$ independent of $n$. Note that $s_n=v_n-w_n$ satisfies
\begin{equation}\label{harmonic}
\begin{cases}
\Delta_{g} s_n=0, & \text { in } B^+_{2r},\\
\dfrac{\partial s_n}{\partial \nu}=0, & \text { on }  \partial B_{2r}^+.
\end{cases}    
\end{equation}

Extending $s_n$ evenly, we have
$$\begin{cases}
\Delta_{g} s_n=0, & \text { in } B_{2r},\\
\dfrac{\partial s_n}{\partial \nu}=0, & \text { on }  \partial B_{2r}.
\end{cases}    $$
Note that for all $y\in B^+_{\frac{r}{2}}(0)$, we see that the reflection of $B_{\frac{r}{2}}(y)\cap  \mathbb R^2_-$ on the horizontal axis is a subset of $B^+_{\frac{r}{2}}(y)$. Here $ \mathbb R^2_-=\{(s,t)\in\mathbb R^2:t<0\}$.
By the mean value property of harmonic functions, this implies that for any $y\in B^+_{\frac{r}{2}}(0)$ we obtain
$$s_n(y)=\fint_{B_{\frac{r}{2}}(y)}s_n(y)dy\leq 2\fint_{B^+_{\frac{r}{2}}(y)}s_n(y)dy,$$
where $\fint$ denotes the mean value. By Jensen's and H\"older inequalities, we have
$$\fint_{B^+_{\frac{r}{2}}(y)} s_n dv\leq\log \fint_{B^+_{\frac{r}{2}}(y)} e^{s_n} dv\leq \frac{1}{2}\log \int_{B^+_{\frac{r}{2}}(y)} e^{2v_n} +\frac{1}{2}\log\int_{B^+_{\frac{r}{2}}(y)} e^{-2w_n} +C\leq C.
$$
Here, $C>0$ is a constant that depends on $B_{\frac{r}{2}}^+(0)$. This gives us a uniform bound $s_n\leq C$ on $B_{\frac{r}{2}}^+(0)$. Therefore, $v_n$ is uniformly bounded from above. Recall that $v_n:=u_n+v_0$. This finishes the proof.
   
\end{proof}

\subsection{Proof of Theorem \ref{theor-blow-up-analysis}}

Given a sequence $\{u_n\}$ of solutions to \eqref{eq017}, we say that a point $p_0\in M$ is a blow-up point for $\{u_n\}$ if for any $r>0$ there holds $\sup_{B_{r}(p)}u_n(x)\rightarrow+\infty$ as $n\rightarrow\infty$. In this case,  Lemma \ref{lem002} implies that $f(p_0)=h(p_0)=0$, since we assume that $f$ and $h$ have the same maximum points. For the sequence $\{u_n\}$ given by Proposition \ref{propo003}, there must exist at least one blow-up point, since otherwise by regularity theory we could extract a subsequence converging smoothly to the global minimizer of $I_{\mu_n,\lambda_n}$, which is a contradiction by \eqref{eq014} and \eqref{eq007}.

As in the proof of Lemma \ref{lemma_bound}, in local coordinates $x$ near $p_0 = 0$, namely $B^+_r:=B^+_r(0)$, for some $r>0,$ there is a smooth function $v_0$ such that $g=e^{2v_0}g_{euc}$.  Thus, if $u_n$ is a solution to \eqref{eq017}, then  $v_n:=u_n+v_0$  satisfies
\begin{equation}\label{eq058}
\begin{cases}
-\Delta v_n=f_{\mu_n} e^{2 v_n} & \text { in }  B_r^+,\\  
\dfrac{\partial v_n}{\partial \nu}=h_{\lambda_n} e^{ v_n} & \mbox{ on }\Gamma_r,   
\end{cases}    
\end{equation}
where $\nu$ is the outward unit normal to $\Gamma_r$. 

\begin{lemma}\label{lemmabounded}
    It holds
    \begin{enumerate}
        \item[(a)] $\displaystyle\limsup_{n\rightarrow\infty}\left(\int_{B^+_r(0)}f_{\mu_n}^+e^{2v_n}dx+\int_{\Gamma_r}h_{\lambda_n}^+e^{v_n}dl\right)\geq \pi$.
        
        \item[(b)] Up to subsequences, the blow-up set for the sequence of solutions $\{u_n\}$ to \eqref{eq017} satisfying \eqref{total_curvature} is finite.
    \end{enumerate}
    \end{lemma}

The proof of item (a) follows by contradiction using \cite[Theorem 1]{MR1132783} and \cite[Lemma 3.2]{MR2504038}, see the proof of \cite[Lemma 2.4]{MR4097234}. Item (b) is a consequence of item (a) and \eqref{total_curvature}.

Let $\{u_n\}$ be the sequence of the solution given by Proposition \ref{propo003}. Let $\{p_\infty^{(i)}:1\leq i\leq N\}$ be the set of blow-up points for $\{u_n\}$. Since $\{u_n\}$ is locally uniformly bounded on $M_{-},$ up to a subsequence, $u_n \rightarrow u_{\infty}$ is smoothly locally on $M_{\infty}:=M\backslash\{p_\infty^{(i)}:1\leq i\leq N\}$. Note that \eqref{total_curvature} gives that  
$\left(-\Delta_{g} u_n\right)_n$ and $\left(\partial u_n/\partial \nu\right)_n$ are  uniformly $L^1$-bounded. Around each $p_\infty^{(i)}$, by \eqref{eq058}, the correspond limit $v_\infty:=u_\infty+v_0$  of $v_n$ satisfies $-\Delta v_\infty=f e^{2 v_\infty}$ in $B_r^+\backslash\{0\}$ and $\partial v_\infty/\partial\nu=he^{v_\infty}$ on $\Gamma_r\backslash\{0\}$.
Using a Green's function representation, as in \eqref{eq037}, we get
$$\lim_{|x|\to 0}\frac{v_\infty(x)}{\log|x|}=-a_i< 0.$$
This implies that $u_\infty(x)=-a_i\log|x|+O(|x|)$ as $x\to 0$. Also, $u_\infty$ satisfies the equation
$$\begin{cases}
    -\Delta_gu_\infty-1=fe^{2u_\infty} & \mbox{ in }M_\infty\\
    \dfrac{\partial u_\infty}{\partial\nu}=he^{u_\infty} & \mbox{ on }\partial M_\infty.
\end{cases}$$
Applying integration by parts twice over $M_r:=M\backslash\cup B_r^+(p_\infty^{(i)})$, using the above equation and then take the limit $r\to 0$ yields
$$\pi\sum_{i=1}^{N}  a_i\phi(p_{\infty}^{(i)})=\int_M(-u_\infty\Delta_{g}\phi-\phi-fe^{2u_\infty}\phi)  +\int_{\partial M}\left( u_\infty\frac{\partial \phi}{\partial \nu}-\phi he^{u_{\infty}}\right),$$
for all $\phi\in C^{\infty}(M).$  
Note that  $a_i \geq 1$  by Lemma \ref{lemmabounded}. Observe that $\{p_\infty^{(i)}:1\leq i\leq N\}\subset \partial M$, since all zeros of $f$ belong to $\partial M$.

\begin{proposition}\label{complete}
 There holds
\begin{itemize}
    \item[(a)] $a_i\in [1,2]$ for every $ 1 \leq i \leq N;$
    \item[(b)]  The metric $g_\infty=e^{2u_\infty}g$ on $M_\infty$ is complete with finite total curvature.
\end{itemize}
\end{proposition}
\begin{proof}
  Consider for each $ 1 \leq i \leq N$ a local conformal chart  $\Psi: B^+_R(p^{(i)}_{\infty}) \rightarrow B^+_R=B_R^+(0) \subset \mathbb{R}_+^2$ with coordinates $x$, such that (locally) $g=e^{2 v_0} g_{euc},$ where $v_0$ is a smooth function defined on $ \overline{B^+_R}$.  If $v_{\infty}(x) = u_{\infty}(x) + v_0(x)$, then $v_{\infty}(x)=-a_i\log|x|+w_{\infty}(x),$ where  $w_{\infty}$ solves
\begin{equation}\label{neweqhelp}
\begin{cases}
-\Delta w_{\infty}=fe^{2v_{\infty}} & \text { in } B^+_R, \\  \dfrac{\partial w_{\infty}}{\partial \nu}=he^{v_{\infty}} & \text { on }  \Gamma_R.\end{cases}
\end{equation}
By \eqref{total_curvature} and Fatou's Lemma, we have $fe^{2v_{\infty}}\in L^1(B_R^+)$ and  $he^{v_{\infty}}\in L^1(\Gamma_R)$. In view of \cite[Theorem 1]{MR1132783} and \cite[Lemma 3.2]{MR2504038}, we can follow the argument from the proof of \cite[Lemma 3.4]{MR3249809} to show that  $e^{2\left|w_{\infty}\right|}$ and $e^{\left|w_{\infty}\right|}$ are $L^p$-bounded in a small neighborhood around $x=0$, for some $p>1$. By the assumption of the maximum point of $f$, we have
$
\alpha^{-1}|x|^2\leq |f(x)|\leq \alpha|x|^2,
$
for some $\alpha>0.$ Therefore
\begin{equation}\label{b_0}
\alpha^{-1}|x|^{2-2a_i} e^{ 2w_{\infty}} \leq\left|f(x)\right| e^{2v_{\infty}} \leq \alpha|x|^{2-2a_i} e^{2w_{\infty}}.
\end{equation}

Fix any $q\in(1,2)$. By Hölder's inequality and \eqref{total_curvature}, as in (5.11) of \cite{MR3351750}, we have the following estimate
\begin{align}
\int_{B^+}|x|^{\frac{2-2a_i}{q}} dx
&=\int_{B^+}(|x|^{2-2a_i}e^{2w_\infty})^{\frac{1}{q}}e^{-\frac{2w_\infty}{q}}dx \nonumber\\
& \leq \alpha^{\frac{1}{q}}\left(\int_{B^+}\left|f(x)\right| e^{2v_{\infty}} d x\right)^{\frac{1}{q}}\left(\int_{B^+} e^{\frac{2|w_\infty|}{q-1}} d x\right)^{\frac{q-1}{q}}\leq C(q) \nonumber,    
\end{align} 
for some small half ball $B^+$ around $x=0$. Thus  $1\leq a_i\leq 2$.

To prove item (b), let us divide the analysis in two cases.

\medskip

\noindent\textbf{Case 1:} $a_i< 2.$

\medskip

For some $q>1$ there holds $\Delta w_\infty\in L^q(B^+),$ where $B^+$ is a small half-ball around $x=0$. Since we also are assuming that  $x=0$ is a non-degenerate maximum point of $h$, we have 
$\beta^{-1}|x|^2\leq |h(x)|\leq \beta|x|^2,$
for some $\beta>0.$ This implies that
\begin{equation}\label{bbb}
\beta^{-1}|x|^{2-a_i} e^{ w_{\infty}} \leq\left|h(x)\right| e^{v_{\infty}} \leq \beta|x|^{2-a_i} e^{w_{\infty}}.
\end{equation}
Since $a_i<2,$ there is a certain small domain $B^+$ around $x=0$ such that, for some $q>1$,   
$\partial w_{\infty}/\partial \nu\in L^q(\Gamma)$ and $\Delta w_\infty\in L^q(B^+).$ By  standard elliptic theory, we have $\|w_{\infty}\|_{L^\infty(\overline{B}^+)}\leq C$ uniformly  on a sufficiently small half-ball $B^+$ around $x = 0$.
 Then there exists a constant $C>0$ such that $e^{2v_\infty} \geq C |x|^{-2a_i} \geq C |x|^{-2}$ near $x=0$ on $B^+$. Consequently, the metric  $g_{\infty} = e^{2u_{\infty}} g = e^{2v_{\infty}} g_{\mathbb{R}^2}$ on $B^+ \setminus \{0\}$ is complete. Since we may also assume pointwise convergence almost everywhere, by Fatou’s Lemma and \eqref{total_curvature}  we get
$$
\int_M\left|f\right| e^{2 u_{\infty}} dv_{g} +\int_{\partial M}\left|h\right| e^{u_{\infty}} da_{g} <\infty.
$$
Then  $g_{\infty}$ has finite total curvature.

\medskip

\noindent\textbf{Case 2:} $a_i=2.$

\medskip

First, similarly to (5.13) in \cite{MR1132783}, using \eqref{neweqhelp} and \eqref{b_0}, for any $d>0$ we get
 $$-\Delta(|x|^{2d} e^{-2w_{\infty}}) \leq  C|x|^{2d-2}.  $$
Also, using \eqref{bbb} and that $\langle x,\nu\rangle=0$ on $\Gamma$, we obtain
$$\frac{\partial}{\partial \nu} \left(|x|^d e^{-w_\infty}\right)= -|x|^d e^{-w_\infty} \frac{\partial w_\infty}{\partial\nu}+d |x|^{d-2}e^{-w_\infty}\langle x,\nu\rangle
\leq C|x|^{d}.$$
Thus, there holds
$$\frac{\partial}{\partial \nu} \left(|x|^{2d} e^{-2w_\infty}\right)=C|x|^{2d}e^{-w_\infty}
\leq C|x|^{2d}e^{|w_\infty|},$$
where the right hand side is in $ L^{q}(\Gamma)$  for some  $q=q(d)>1$, provided  on a sufficiently small half-ball around $x = 0$ we have $e^{|w_\infty|}\in L^q.$  Since we also have $\Delta(|x|^{2d} e^{-2w_{\infty}})\in L^p(B^+)$ for some $p=p(d)>1$ we have by elliptic estimates that  $|x|^d e^{-w_{\infty}} \leq C$.  Thus, for any $k$ there is a constant $k>0$ such that $e^{2v_{\infty}}=|x|^{-4} e^{2w_{\infty}} \geq k|x|^{d-4}$, and hence  $g_{\infty}$  is complete on $B^+ \backslash\{0\}$.
\end{proof}

Motivated by \cite[Proposition 5.1]{MR4517687}, in the proof of Theorem \ref{theor-blow-up-analysis} we will use Ekeland’s variational principle, which we recall below for the sake of the reader. See \cite[Theorem 5.1]{MR2431434} pg. 51 and \cite[Theorem 1.4.1]{MR1921556}.

\begin{theorem}\label{teo-ekeland}
    Let $(X,d)$ be a complete metric space and consider a function $\varphi:X\to(-\infty,+\infty]$ that is lower semi-continuous, bounded from below and not identical to $+\infty$. Let $\varepsilon>0$ and $\lambda>0$ be given and let $x\in X$ be such that $\varphi(x)\leq \inf_X\varphi+\varepsilon$. Then there exists $x_\varepsilon\in X$ such that
    \begin{enumerate}[(a)]
        \item $\varphi(x_\varepsilon)\leq\varphi(x)$,
        \item $d(x_\varepsilon,x)\leq\lambda$,
        \item $\varphi(x_\varepsilon)<\varphi(z)+\frac{\varepsilon}{\lambda}d(x_\varepsilon,z)$ for every $z\not=x_\varepsilon$.
    \end{enumerate}
\end{theorem}

\begin{proof}[Proof of Theorem \ref{theor-blow-up-analysis}]  To conclude our analysis, our goal is to analyze the behavior of the solution as we approach each point $p_{\infty}^{(i)}$, for $1\leq i \leq N$. We will consider two sequences $\lambda_n\downarrow 0$ and $\mu_n\downarrow 0$, with $\lambda_n^2-\lambda_n^3\leq\mu_n\leq \lambda_n^2+\lambda_n^3$, such that $c_{\mu_n,\lambda_n}'(0)$ exists, and the corresponding sequence of solutions $u_n$, which are given by  Proposition \ref{propo003}. As illustrated in the proof of Proposition \ref{complete}, consider a local conformal chart $\Psi: B^+_R(p^{(i)}_{\infty}) \rightarrow B^+_R:=B^+_R(0) \subset \mathbb{R}_+^2$ with coordinates $x=(s,t)$ and $\Psi(p_\infty^{(i)})=0$, and let $v_{n}(x) = u_{n}(x) + v_0(x)$ be such that (locally) $g=e^{2 v_0} g_{euc},$ where $v_0$ is a smooth function on $ \overline{B^+_R}$. Here $R$ is chosen sufficiently small in order to guarantee that in $B_R^+(p_\infty^{(i)})$ the only maximum point of $f$ and $h$ is $p_\infty^{(i)}$. The function $v_n$ solves
\begin{equation}\label{eq027}
    \begin{cases}
-\Delta v_n=f_{\mu_n}(x) e^{2 v_n} & \text { in }  B_R^+,\\  
\dfrac{\partial v_n}{\partial \nu}=h_{\lambda_n} e^{ v_n} & \mbox{ on }\Gamma_R.    
\end{cases}
\end{equation}
 Since $p_\infty^{(i)}$ is a non-degenerate  maximum point of $f$ and $h$, we have
\begin{equation}\label{eq028}
    f(x)=\frac{1}{2}\nabla^2f(p_\infty^{(i)})(x,x)+O(|x|^3),\;\mbox{ for all } x\in B_R^+
\end{equation}
and
\begin{equation}\label{eq032}
h(x)=\frac{1}{2}\nabla^2h(p_\infty^{(i)})(x,x)+O(|x|^3),\;\mbox{ for all } x\in\Gamma_R,
\end{equation}
with $\nabla^2f(p_\infty^{(i)})$ and $\nabla^2h(p_\infty^{(i)})$ negative definite. We can choose $R$ small enough such that
\begin{equation}\label{eq030}
    -\alpha_1|x|^2\leq f(x)\leq -\alpha_2|x|^2,\;\mbox{ for all } x\in B_R^+(0)
\end{equation}
and
\begin{equation}\label{eq054}
    -\alpha_3|x|^2\leq h(x)\leq -\alpha_4|x|^2,\;\mbox{ for all } x\in \Gamma_R,
\end{equation}
with $0<\alpha_2\leq\alpha_1$ and $0<\alpha_4\leq\alpha_3$. Define
$K_n:=\{p\in M:f(p)+\mu_n\geq 0\}\cap \overline{B_R^+(p_\infty^{(i)})}$
and $Q_n:=\{p\in\partial M:h(p)+\lambda_n\geq 0\}\cap \overline{B_R^+(p_\infty^{(i)})}.$
Since $f(p_\infty^{(i)})=h(p_\infty^{(i)})=0$, then $K_n$ and $Q_n$ are nonempty sets. Now, for each $n$ consider $a_n=\sqrt{\mu_n/\alpha_1}$. Using that $\mu_n\leq \lambda_n^2+\lambda_n^3$ and the inequalities \eqref{eq030} and \eqref{eq054} we obtain $B_{a_n}^+:=B_{a_n}^+(p_\infty^{(i)})\subset K_n$ and $\Gamma_{a_n}\subset Q_n$, for $n$ large enough. Consider a sequence $(x_{n})$ in $M$ such that $v_{n}(x_{n}) = \sup_{B_{a_n}^+}v_{n}(x)$. Since $f(x_n)+\mu_n\geq 0$, by \eqref{eq030} there exists some constant $C>0$ independent of $n$ such that
\begin{equation}\label{eq029}
    |x_n|^2\leq C\mu_n.
\end{equation}
Up to a subsequence, we see that we have two mutually disjoint cases:

\medskip

    \noindent{\bf Case 1:} $\lambda_{n}^{2}e^{v_{n}(x_{n})} \rightarrow +\infty$.

\medskip

In this case, for $x \in B_{a_n}^{+},$ rescale
\begin{equation*}
w_{n}(x) = v_{n}(x_n + r_{n}x) - v_{n}(x_{n})\leq w_n(0)=0,
\end{equation*}
which is defined in $\Omega_n:=\{x\in \mathbb R^2:x_n+r_nx\in B_{a_n}^+(0)\}$, where $r_n>0$ is so that 
$$r_{n}\lambda_{n}e^{v_{n}(x_{n})} = 1.$$
 Note that $r_n/\lambda_n\rightarrow 0$, and using that $\lambda_n^2-\lambda_n^3\leq\mu_n$ we have $a_n/r_n\rightarrow+\infty$. Passing to a subsequence, we can assume that
 $\frac{d(x_n,\Gamma_n)}{r_n}\rightarrow t_0,$
 with either $t_0\geq 0$ or $t_0=+\infty$. 
Since $v_n$ satisfies \eqref{eq027}, we have
$$
\begin{cases}-\Delta w_n = r_n^2f_{\mu_n}(x_n+r_nx)e^{2w_n}e^{2v_n(x_n)}=F_ne^{2w_n}& \text { in } \Omega_n, \\ 
\dfrac{\partial w_n}{\partial \nu}  = r_nh_{\lambda_n}(x_n+r_nx)e^{w_n}e^{v_n(x_n)}=H_ne^{w_n} & \text { on } \Gamma_n,\end{cases}
$$
where $\Gamma_n$ is the straight portion of $\partial \Omega_n$, $F_n(x):=f_{\mu_n}(x_n+r_nx)/\lambda_n^2$ and $H_n(x):=h_{\lambda_n}(x_n+r_nx)/\lambda_n$. Set $A=\frac{1}{2}\nabla^2f(p_\infty^{(i)})$ and $B=\frac{1}{2}\nabla^2h(p_\infty^{(i)})$. Using \eqref{eq028} and \eqref{eq032},  we obtain
$$f(x_n+r_nx)-f(x_n) =2r_nA(x_n,x_n)+r_n^2A(x,x)+O(|x_n|^3)+O(|r_nx|^3)$$
and
$$h(x_n+r_nx)-h(\tilde x_n) =2r_nB(\tilde x_n,\tilde x_n)+r_n^2B(s,s)+O(|\tilde x_n|^3)+O(|r_ns|^3).$$
Here $x=(s,t)$ and $\tilde x_n=(x_{n1},0)$, where $x_n=(x_{n1},x_{n2})$. By \eqref{eq029}, for $|x|\leq C$,  we obtain
$$\frac{|f(x_n+r_nx)-f(x_n)|}{\mu_n}\leq C\left(r_n+\frac{r_n^2}{\mu_n}+|x_n|\right)$$
and
$$\frac{|h(x_n+r_nx)-h(\tilde x_n)|}{\lambda_n}\leq C\left(\frac{r_n}{\lambda_n}+\frac{\mu_n}{\lambda_n}|\tilde x_n|\right).$$
Since $\mu_n/\lambda_n^2\to 1$, we get
$$\lim_{n\rightarrow+\infty}F_n(x)=\lim_{n\rightarrow+\infty}\frac{\mu_n}{\lambda_n^2}\frac{f_{\mu_n}(x_n+r_nx)}{\mu_n}=\lim_{n\rightarrow+\infty}\frac{f(x_n)}{\mu_n}+1
$$
and
$$\lim_{n\rightarrow+\infty}H_n(x)=\lim_{n\rightarrow+\infty}\frac{h(\tilde x_n)}{\lambda_n}+1.$$
Note that $f\leq 0$ and $h\leq 0$ imply that both limits are less than equal to 1. On the other hand, $f_{\mu_n}(x_n+r_nx)\geq 0$, for $x\in\Omega_n$, and $h_{\lambda_n}(x_n+r_nx)\geq 0$, for $x\in\Gamma_n$, which implies that both limits are nonnegative. Therefore, up to subsequences, both limits exist. Let $\lim_{n\to\infty} F_n(x)=c_\infty$ and $\lim_{n\to\infty} H_n(x)=d_\infty$, for some $c_\infty,d_\infty\in[0,1]$. Besides
\begin{equation}\label{eq049}
    \int_{\Omega_n}e^{2w_n(x)}dx \leq r_n^{-2}e^{-2v_n(x_n)}\int_Me^{2v_n}dv_g=\frac{\lambda_n^2}{\mu_n}\mu_n\int_Me^{2v_n}dv_g
\end{equation}
and
\begin{equation}\label{eq050}
\int_{\Gamma_n}e^{w_n(x)}dl = \int_{\Gamma_n}e^{v_n(x_n+r_nx)-v_n(x_n)}dl\leq \lambda_n\int_{\partial M}e^{v_n}d\sigma_g.    
\end{equation}
Since $\lambda_n^2/\mu_n\rightarrow 1$, it follows by \eqref{eqbounds} that
$$\frac{1}{2}\int_{\Omega_n}e^{2w_n}dx+\int_{\Gamma_n}e^{w_n}dl\leq 2\pi.$$

If $t_0=+\infty$, then $\Omega_n$ exhausts $\mathbb R^2$. It follows from the proof of Theorem 1.4 of \cite{MR3351750} that $w_n\rightarrow w_\infty$ in $C^2_{loc}(\mathbb R^2)$, where $w_\infty$ satisfies the Liouville equation 
$$-\Delta w_\infty=c_\infty e^{2w_\infty}\quad\mbox{ in }\mathbb R^2,
\mbox{ with }\int_{\mathbb R^2}e^{2w_\infty}dx\leq 4\pi.$$
This implies that $c_\infty>0$. Considering $r_n/\sqrt{c_\infty}$ instead of $r_n$, we obtain that $w_\infty$ solves
 $$-\Delta w_\infty= e^{2w_\infty}\quad\mbox{ in }\mathbb R^2,
\mbox{ with }\int_{\mathbb R^2}e^{2w_\infty}dx\leq 4\pi.$$
From the classification result in \cite{MR1121147}, it follows that $w_\infty(x)=-\log\left(1+\frac{1}{8}|x|^2\right)$ with
$$  \int_{\mathbb R^2}e^{2w_\infty}dx=4\pi.  $$
Thus we have
$$\int_{\Omega_n}F_ne^{2w_n}dx=4\pi+o(1),$$
as $n\rightarrow+\infty$. On the other hand, after a change of variables, \cite[Lema 3.7]{MR4546499} implies that 
$$\int_{\Omega_n}F_ne^{2w_n}dx+\int_{\Gamma_n}H_ne^{w_n}dl\leq 2\pi+o(1),$$
which is a contradiction. Therefore, $t_0<+\infty$. In this case, $\Omega_n$ exhausts $\mathbb R_{t_0}^2:=\mathbb R\times(t_0,+\infty)$. By standard elliptic theory, see \cite{MR4553958}, for instance, up to a subsequence, $w_n\rightarrow w_\infty$ in $C_{loc}^2(\mathbb R_{t_0}^2)$, where $w_\infty$ satisfies
\begin{equation}\label{eq035}
    \begin{cases}
    -\Delta w_\infty= c_\infty e^{2w_\infty} & \mbox{ in }\mathbb R^2_{t_0},\\
    \dfrac{\partial w_\infty}{\partial\nu}=d_\infty e^{w_\infty} & \mbox{ on }\partial\mathbb R^2_{t_0},
\end{cases}\;\;\mbox{ with }\;\;\frac{1}{2}\int_{\mathbb R^2_{t_0}}e^{2w_\infty}dx+\int_{\partial \mathbb R^2_{t_0}}e^{w_\infty}dl\leq 2\pi,
\end{equation}
and $c_\infty,d_\infty\in[0,1]$. Thanks to the classifications of solutions to the equation \eqref{eq035}, see \cite{MR1369398,MR1971036}, we obtain the existence of $\Lambda>0$ and $s_0 \in \mathbb{R}$ such that
    $$w_\infty(s,t)=\log \frac{2\Lambda}{c_\infty\Lambda^2+\left(s-s_0\right)^2+\left(t-t_0+d_\infty\Lambda\right)^2},$$
where $d_\infty>0$ if $c_\infty=0$. This finishes the analysis of blow-up in Case 1. 

Finally, we prove an estimate for the number of blow-up points. Before, recall that

$$
d_\infty \int_{\partial \mathbb{R}_{+}^2} e^{w_{\infty}}dl=\beta, \quad c_\infty \int_{\mathbb{R}_{+}^2} e^{2w_{\infty}}dx=2 \pi-\beta,
$$
where $\beta:=2 \pi \frac{d_\infty}{\sqrt{d_\infty^2+c_\infty}}\leq 2\pi$ (see for instance Section 2.2 of \cite{MR4546499}).
As in \eqref{eq049} and \eqref{eq050}, at each blow-up point $p_{\infty}^{(i)},$  we must have
$$2\pi-\beta= c_{\infty}\int_{\mathbb{R}_+^2} e^{2 w_{\infty}} d x   \leq c_{\infty} \limsup _{n \rightarrow \infty} \int_{\Omega_n} e^{2 w_n} d x   \leq c_{\infty} \limsup _{n \rightarrow \infty}\left( \frac{\lambda_n^2}{\mu_n}\mu_n \int_{B_{a_n}^+(p_\infty^{(i)})} e^{2 u_n} dv_g\right)
$$
and
$$
\beta\leq d_{\infty} \limsup _{n \rightarrow \infty}\left(\mu_n\int_{\Gamma_{a_n}(p_\infty^{(i)})} e^{ u_n} d\sigma_g\right).
$$
Since $p_\infty^{(i)}$ is the only maximum point in  $B_R^+(p_\infty^{(i)})$, we have
$$(2\pi-\beta)N\leq c_{\infty} \limsup _{n \rightarrow \infty}\left( \frac{\lambda_n^2}{\mu_n}\mu_n \int_{M} e^{2 u_n} dv_g\right)\quad\mbox{ and }\quad\beta N\leq d_{\infty} \limsup _{n \rightarrow \infty}\left(\mu_n\int_{\partial M} e^{ u_n} d\sigma_g\right),$$
where $N\geq 1$ is the number of blow-up point. Consider the following cases
\begin{itemize}
    \item If $c_\infty=0$, then $\beta=2\pi$ and $d_\infty=1$. 
In this case, $N=1$. In fact, by \eqref{eqbounds}, it holds
    $$
2\pi N\leq \limsup _{n \rightarrow \infty}\left(\mu_n \int_{\partial M} e^{u_n} d\sigma_g\right)\leq 2\pi.
$$
\item  If $d_\infty=0$, then $\beta=0$ and  $c_\infty\in (0,1]$. In this case, $N\leq 2$. In fact, again by \eqref{eqbounds}, it holds
$$
      2\pi N\leq c_{\infty} \limsup _{k \rightarrow \infty}\left( \frac{\lambda_n^2}{\mu_n}\mu_n \int_{M} e^{2 u_n} dv_g\right)\leq  4\pi c_{\infty}.
$$

      \item If $c_\infty, d_\infty\in (0,1],$ then $N\leq 3$. In fact, as before
      \begin{align*}
2\pi N\leq c_{\infty} \limsup _{k \rightarrow \infty}\left( \frac{\lambda_n^2}{\mu_n}\mu_n \int_{M} e^{2 u_k} d v_g\right)+d_{\infty} \limsup _{n \rightarrow \infty}\left(\lambda_n \int_{\partial M} e^{u_k} d \sigma_g\right)\leq 6\pi. 
      \end{align*}

\end{itemize}

\medskip

    \noindent{\bf Case 2:} $\lambda_{n}^{2}e^{v_{n}(x_{n})} \leq C$.

\medskip


In view of \eqref{eq030} we obtain $f(x)+\mu_n\leq 0$ if $\mu_n\leq c|x|^2$, for some positive constant $c$. Thus, by Lemma \ref{lemmabounded} we obtain
$$\pi\leq\int_{B^+_r(0)}f_{\mu_n}^+(x)e^{2v_n(x)}dx=\int_{B^+_{c\sqrt{\mu_n}}(0)}f_{\mu_n}(x)e^{2v_n(x)}dx\leq c\mu_n^{2}e^{2v_n(x_n)},$$ 
which implies that $c\mu_ne^{v_n(x_n)}\geq 1$, for some constant $c>0$ independent of $n$. Since $\lambda_n^2/\mu_n\rightarrow 1$, then $    ce^{-v_n(x_n)}< \lambda_n^2< Ce^{-v_n(x_n)},$ for positive constants $c$ and $C$, which does not depend on $n$. By continuity, for each $n$ we can find $\varepsilon_n\in(0,\lambda_n^2)$ such that 
\begin{equation}\label{eq055}
ce^{-v_n(z)}\leq \lambda_n^2\leq Ce^{-v_n(z)},\quad\mbox{for all } z\in B_R^+(0)\mbox{ with } |z-x_n|<\varepsilon_n
\end{equation}
By Ekeland’s variational principle (Theorem \ref{teo-ekeland}), taking $\varphi=e^{-\frac{v_n}{2}}$, $\lambda=\varepsilon_n$ and $\varepsilon=\varepsilon_n^2$, we can find a sequence $y_n\in B_R^+$ such that
\begin{enumerate}[(a)]
    \item $v_n(x_n)\leq v_n(y_n)$,
    \item $|x_n-y_n|<\varepsilon_n$,
    \item $e^{-v_n(y_n)}<e^{-v_n(x)}+\varepsilon_n|x-y_n|$, for every $x\not=y_n$.
\end{enumerate}

By $(a)$ and $(b)$, it holds $v_n(y_n)\to+\infty$ and $y_n\rightarrow 0$. Let $w_n(y)=v_n(\lambda_ny)+2\log\lambda_n$, defined in $B_n:=B_{R/\lambda_n}^+(0)$. Note that $w_n$ satisfies
$$
\begin{cases}-\Delta w_n =F_ne^{2w_n}& \text { in } B_n, \\ 
\dfrac{\partial w_n}{\partial \nu}  = H_ne^{w_n} & \text { on } \Gamma_n,
\end{cases}
$$
where $F_n(y)=f(\lambda_ny)/\lambda_n^2+\mu_n/\lambda_n^2$ and $H_n(y)=h(\lambda_ny)/\lambda_n+1$. Since $\mu_n/\lambda_n^2\rightarrow 1$,  by \eqref{eq028} and \eqref{eq032} we obtain that $F_n(y)\rightarrow\frac{1}{2}\nabla^2f(p_\infty^{(i)})(y,y)+1$ and $H_n(y)\rightarrow1$. 

By \eqref{eq029}, we have $|x_n/\lambda_n|<L$, for some $L>0$. Note that if $|x|<\lambda_nr$, then $|x-y_n|\leq \lambda_n L_0$, with $L_0:=r+L+1$. 

\medskip\noindent{\bf Claim:} Given $r>0$, it holds $w_n(y)<c$,  for all $y\in B_n$ with $|y|\leq r$, for $n$ large enough, where $c>0$ does not depend on $n$.

\medskip

Using \eqref{eq055}, item (c) and that $|x_n/\lambda_n|<L$, if $|x|\leq \lambda_nr$, we get
$$\lambda_n^2e^{v_n(x)}\leq C(1-CL_0\lambda_n)^{-1},$$
which implies the claim.

\medskip

By this claim and the Harnack-type inequalities (see \cite[Lemma A.2]{MR3249809}), $w_n$ is uniformly bounded in $L_{loc}^\infty(\mathbb R^2_+)$. Therefore, up to a subsequence,
$$w_n\rightarrow w_\infty\quad\mbox{ in }C^2_{loc}(\mathbb R^2_+),$$
which is a solution of the equation
$$
\begin{cases}
-\Delta w_\infty =(1+\frac{1}{2}\nabla^2f(p_\infty^{(i)})(y,y))e^{2w_\infty} & \text { in } \mathbb R^2_+, \\ 
\dfrac{\partial w_\infty}{\partial \nu}  = e^{w_\infty} & \text { on } \partial \mathbb R^2_+.
\end{cases}
$$

As before and using \eqref{total_curvature} we have
$$\int_{\mathbb R^2_+}e^{2w_\infty}dx+\int_{\partial R_+^2}e^{w_\infty}dl<\infty\quad\mbox{ and }\quad \int_{\mathbb{R}^{2}_{+}}\left|1+ \frac{1}{2}\nabla^2f(p_\infty^{(i)})(y,y)\right|e^{2w_\infty}dx<\infty.$$
By Theorem \ref{teo001}, whose proof we postponed to the next section, we conclude that this is impossible. Therefore, we complete the proof of Theorem \ref{theor-blow-up-analysis}.
\end{proof}

\subsection{Ruling out slow blow-up: Nonexistence result}\label{rule}

In this subsection, our primary objective is to eliminate certain blow-up cases identified in the latter part. We aim to achieve this by conducting an analysis akin to the one presented in \cite{MR4753069} and \cite{MR4153107}, specifically excluding cases where the blow-up occurs in a ``slow" manner.

\begin{theorem}[Theorem \ref{teo001}]\label{ruling}
Suppose that $A$ is a negative definite $2\times 2$ matrix.  Then there is no solution $w \in C^{\infty}(\mathbb{R}^{2}_{+})$ of the equation
\begin{equation}\label{eq0000}
\begin{cases}
    -\Delta w=Fe^{2w} & \mbox{ in }\R^{2}_{+},\\
    \dfrac{\partial w}{\partial\nu}=e^w & \mbox{ on }\partial \mathbb R^{2}_{+},
\end{cases}
\end{equation}
with $w\leq C$, where $F(x):=1+ (Ax,x)$, such that $Fe^{2w}\in L^1(\mathbb R^2_+)$,
\begin{equation*}
   V_0:= \int_{\R^{2}_{+}}e^{2w} dx<\infty,\quad H_0:= \int_{\partial\R^{2}_{+}}e^{w} dl<\infty,\quad \mbox{ and }\quad    K_{0} := \int_{\mathbb{R}^{2}_{+}}Fe^{2w} dx\in \R.
\end{equation*}
\end{theorem}

\begin{proof}  By contradiction, suppose $w$ is a solution to \eqref{eq0000}. Motivated by the proof of \cite[Lemma 2.1]{MR4753069} and \cite[Theorem 4.5]{MR4153107}, consider the function
\begin{align}\label{auxiliar}
    \tilde{w}(x) =& -\frac{1}{2\pi}\int_{\mathbb{R}^{2}_{+}}(\log|x-y| + \log|\overline{x}-y|- 2\log|y|)Fe^{2w}dy   \\ 
    &  - \frac{1}{2\pi}\int_{\partial \mathbb{R}^{2}_{+}}(\log|x-y| + \log|\overline{x}-y| -2\log|y|)e^w dy,\nonumber
\end{align}
where $\overline{x}$ is the reflection point of $x$ about the boundary $\partial\R^{2}_{+}$, that is, if $x=(s,t)$, then $\overline x=(s,-t)$. As in \cite[Proposition 3.2]{MR1369398}, it holds
\begin{equation}\label{eq037}
    \lim_{|x|\rightarrow \infty}\frac{\tilde{w}(x)}{\log|x|} = -\frac{1}{\pi}\int_{\mathbb{R}^{2}_{+}}Fe^{2w} - \frac{1}{\pi}\int_{\partial \mathbb{R}^{2}_{+}}e^w
= -\frac{1}{\pi}(K_{0} + H_{0})=:-d.
\end{equation}
It is easy to check that $\tilde{w}$ is a solution of 
$$
\begin{cases}
    -\Delta \tilde w=Fe^{2 w} & \mbox{ in }\R^{2}_{+},\\
    \dfrac{\partial \tilde  w}{\partial\nu}=e^w & \mbox{ on }\partial \R^{2}_{+}.
\end{cases}
$$
In this way, the function $v:= w - \tilde{w}$ is harmonic in $\mathbb{R}^{2}_{+}$ with $\partial v/\partial \nu = 0$ in $\partial \mathbb{R}^{2}_{+}$. We extend $v$ to $\mathbb R^2$ by even reflection. Since $|v|\leq C+C\log(2+|x|)$ for some constant $C$ (see \cite[Proposition 3.1]{MR1369398}), we find that $v$ is constant. Hence $w=\tilde w + C$ for some constant $C$.  

From \eqref{eq037}, for all $\varepsilon>0$ we have $e^{\tilde w}\geq |x|^{-d-\varepsilon}$, for all $|x|\gg 1$. Since $A$ is a negative definite matrix, we obtain
\begin{align*}
 \int_{\mathbb{R}_+^2 \backslash B^+_1}|x|^{2-2 d-2\varepsilon} d x &\leq c \int_{\mathbb{R}_+^2}|x|^2 e^{2 \tilde{w}} d x+c  \leq c \int_{\mathbb{R}_+^2}|(A x, x)| e^{2 w} d x+c\\
 &=-c \int_{\mathbb{R}_+^2}(A x, x) e^{2 w} d x+c\\
 & =c \int_{\mathbb{R}_+^2} e^{2 w} d x-c \int_{\mathbb{R}_+^2}F e^{2 w} d x=c\left(V_0-K_0\right)+c.    
\end{align*}
Thus $d \geq 2$ and consequently $K_0+H_0 \geq  2 \pi$. 

\medskip

\noindent{\bf Claim 1:} There exists a positive constant $C>0$ such that for $|x|\gg 1$ it holds
$$\tilde w(x)\leq -d\log|x|+C.$$

First, using \eqref{auxiliar} and \eqref{eq037}, a direct computation yields that
\begin{align*}
    \tilde w(x)&+d\log|x| = \frac{1}{2\pi}\int_{\mathbb R^2_+}\left(\log\frac{|x|(|y|+1)}{|x-y|}+\log\frac{|x|(|y|+1)}{|\overline x-y|}\right)F^+e^{2w}dy\\
    & -\frac{1}{2\pi}\int_{\mathbb R^2_+}\left(\log\frac{|x|(|y|+1)}{|x-y|}+\log\frac{|x|(|y|+1)}{|\overline x-y|}\right)F^-e^{2w}dy\\
    &+ \frac{1}{\pi}\int_{\partial\mathbb R^2_+}\log\frac{|x|(|y|+1)}{|x-y|}e^{w}dl+ \frac{1}{\pi}\int_{\mathbb R^2_+}\log\frac{|y|}{|y|+1}Fe^{2w}dy+ \frac{1}{\pi}\int_{\partial\mathbb R^2_+}\log\frac{|y|}{|y|+1}e^{w}dl\\
    & = J_1+J_2+J_3+J_4+J_5.
\end{align*}
In the third term on the right hand side, we used $|x-y|=|\overline x-y|$, for all $y\in\partial\mathbb R^2_+$. Here $F^{\pm}=\max\{\pm F,0\}\geq 0$. Note that for $|x|\geq 1$, we have $|x-y|\leq |x|+|y||x|$ and $|\overline x-y|\leq |x|+|y||x|$. This implies $J_2\leq 0$. Also, $J_5\leq 0$, since $\log\frac{r}{r+1}<0$.

For the term $J_4$, note that there exists $\delta>0$ such that $F(y)>c>0$, for all $|y|\leq \delta$. Since $\log\frac{r}{r+1}<0$ and $[\delta,\infty)\ni r\mapsto\log\frac{r}{r+1}$ is a bounded function, we obtain
$$J_4\leq \frac{1}{\pi}\int_{\mathbb R^2_+\backslash B_\delta^+(0)}\log\frac{|y|}{|y|+1}Fe^{2w}dy\leq C\|Fe^{2w}\|_{L^1(\mathbb R^2_+)}<\infty.$$

For the term $J_1$, we split it as
\begin{align*}
    J_1 & =  \frac{1}{2\pi}\int_{B^+_{\frac{|x|}{2}}(x)}\left(\log\frac{|x|(|y|+1)}{|x-y|}+\log\frac{|x|(|y|+1)}{|\overline x-y|}\right)F^+e^{2w}dy\\
    &+  \frac{1}{2\pi}\int_{\mathbb R^2_+\backslash B^+_{\frac{|x|}{2}}(x)}\left(\log\frac{|x|(|y|+1)}{|x-y|}+\log\frac{|x|(|y|+1)}{|\overline x-y|}\right)F^+e^{2w}dy\\
    & = J_{11}+J_{12}
\end{align*}

Note that $F(y)=1+(Ay,y)<1$ for all $y\in\mathbb R^2_+$. From \eqref{eq037}, for all $\varepsilon>0$ we have 
\begin{equation}\label{eq056}
e^{\tilde w}\leq |x|^{-d+\varepsilon},    \quad\mbox{for all } |x|\gg 1.
\end{equation}
In addition, for $|x-y|\leq|x|/2$ we have $|x|/2\leq|y|\leq 3|x|/2$. Then, for all $|x|\gg 1$, yields
\begin{align*}
    J_{11} & =\frac{1}{\pi}\int_{B^+_{\frac{|x|}{2}}(x)}\log(|x|(|y|+1))F^+e^{2w}dy\\
    &+\frac{1}{2\pi}\int_{B^+_{\frac{|x|}{2}}(x)}\left(\log\frac{1}{|x-y|}+\log\frac{1}{|\overline x-y|}\right)F^+e^{2w}dy\\
    & \leq C|x|^{-2d+2\varepsilon}\log\left(\frac{3}{2}|x|^2+|x|\right)\int_{B^+_{\frac{|x|}{2}}(x)}dy\\
    & +C|x|^{-2d+2\varepsilon}\int_{B^+_{\frac{|x|}{2}}(x)}\left(\frac{1}{|x-y|}+\frac{1}{|\overline x-y|}\right)dy\\
    & \leq C|x|^{2-2d+2\varepsilon}\log\left(\frac{3}{2}|x|^2+|x|\right)+C|x|^{1-2d+2\varepsilon}\leq C,
\end{align*}
since $d\geq 2$. In the last inequality, we used $|\overline x-y|\geq|x-y|$, for all $x,y\in\mathbb R^2_+$. For $J_{12}$, note that $|x|/|x-y|\leq 2$ if $y\in\mathbb R^2_+\backslash B^+_{\frac{|x|}{2}}(x)$. Thus, using \eqref{eq056} we get
$$J_{12}\leq \frac{1}{\pi}\int_{\mathbb R^2_+\backslash B^+_{\frac{|x|}{2}}(x)}\log\frac{|x|(|y|+1)}{|x-y|}F^+e^{2w}dy\leq C\int_{\mathbb R^2_+\backslash B^+_{\frac{|x|}{2}}(x)}|y|^{-2d+2\varepsilon}\log(2|y|+2)dy\leq C,$$
since $d\geq 2$. Therefore, $J_1\leq C$.

Similarly, we obtain an estimate for the term $J_3$. We divide it as follows.
\begin{align*}
    J_3 & = \frac{1}{\pi}\int_{\Gamma_{\frac{|x|}{2}}(x)}\log\frac{|x|(|y|+1)}{|x-y|}e^wdl+\frac{1}{\pi}\int_{\partial\mathbb R^2_+\backslash\Gamma_{\frac{|x|}{2}}(x)}\log\frac{|x|(|y|+1)}{|x-y|}e^wdl\\
    & =J_{31}+J_{32}.
\end{align*}
As we estimated $J_{11}$ and $J_{12}$, we obtain
\begin{align*}
    J_{31} & \leq C|x|^{1-d+\varepsilon}\log\left(\frac{3}{2}|x|^2+|x|\right)+C|x|^{-d+\varepsilon}\int_{\Gamma_{\frac{|x|}{2}}(x)}\frac{dl}{|x-y|}\\
    &\leq C|x|^{1-d+\varepsilon}\log\left(\frac{3}{2}|x|^2+|x|\right)+C|x|^{-d+\varepsilon}\log|x|\leq C
\end{align*}
and
\begin{align*}
    J_{32} & \leq C\int_{\partial\mathbb R^2_+\backslash\Gamma_{\frac{|x|}{2}}(x)}|y|^{-d+\varepsilon}\log(2|y|+2)dl\leq C.
\end{align*}
This implies Claim 1.

\medskip

Since $\nabla\tilde w=\nabla w$ and $|\overline x-y|=|x-\overline y|$, we have

\begin{align*}
    \langle x,\nabla \tilde w\rangle = & -\frac{1}{2\pi}\int_{\mathbb R^2_+}\left(\frac{\langle x,x-y\rangle}{|x-y|^2}+\frac{\langle x,x-\overline y\rangle}{|\overline x-y|^2}\right)F(y)e^{2w(y)}dy\\
    & -\frac{1}{2\pi}\int_{\partial \mathbb R^2_+}\left(\frac{\langle x,x-y\rangle}{|x-y|^2}+\frac{\langle x,x-\overline y\rangle}{|\overline x-y|^2}\right)e^{w(y)}dl.
\end{align*}
Thus
$$\int_{B_R^+(0)} F(x)e^{2w(x)} \langle x,\nabla w\rangle dx+\int_{\Gamma_R}\langle x,\nabla w\rangle e^wdl$$
$$=    \int_{B_R^+(0)} F(x)e^{2w(x)}\left(-\frac{1}{2\pi}\int_{\mathbb R^2_+}\left(\frac{\langle x,x-y\rangle}{|x-y|^2}+\frac{\langle x,x-\overline y\rangle}{|\overline x-y|^2}\right)F(y)e^{2w(y)}dy\right.$$
$$\left.-\frac{1}{2\pi}\int_{\partial \mathbb R^2_+}\left(\frac{\langle x,x-y\rangle}{|x-y|^2}+\frac{\langle x,x-\overline y\rangle}{|\overline x-y|^2}\right)e^{w(y)}dl\right)dx$$
\begin{align*}
    +\int_{\Gamma_R} e^{w(x)}\left(-\frac{1}{2\pi}\int_{\mathbb R^2_+}\left(\frac{\langle x,x-y\rangle}{|x-y|^2}+\frac{\langle x,x-\overline y\rangle}{|\overline x-y|^2}\right)F(y)e^{2w(y)}dy\right.\\
    \left.-\frac{1}{2\pi}\int_{\partial \mathbb R^2_+}\left(\frac{\langle x,x-y\rangle}{|x-y|^2}+\frac{\langle x,x-\overline y\rangle}{|\overline x-y|^2}\right)e^{w(y)}dl\right)dl.
\end{align*}

For the right hand side of the above identity, expressing $x$ as $x=\frac{1}{2}((x+y)+(x-y))$ and $x=\frac{1}{2}((x+\overline y)+(x-\overline y))$, and taking into account that $|x-\overline y|=|\overline x-y|$, we get
\begin{align*}
    RHS = & \int_{B_R^+(0)} F(x)e^{2w(x)}\left(-\frac{1}{2\pi}\int_{\mathbb R^2_+}F(y)e^{2w(y)}dy-\frac{1}{2\pi}\int_{\partial \mathbb R^2_+}e^{w(y)}dl\right)dx\\
   &+\frac{1}{2} \int_{B_R^+(0)} F(x)e^{2w(x)}\left(-\frac{1}{2\pi}\int_{\mathbb R^2_+}\left(\frac{\langle x+y,x-y\rangle}{|x-y|^2}+\frac{\langle x+\overline y,x-\overline y\rangle}{|\overline x-y|^2}\right)F(y)e^{2 w(y)}dy\right.\\
   & \left.-\frac{1}{2\pi}\int_{\partial \mathbb R^2_+}\left(\frac{\langle x+y,x-y\rangle}{|x-y|^2}+\frac{\langle x+\overline y,x-\overline y\rangle}{|\overline x-y|^2}\right)e^{w(y)}dl\right)dx\\
    &+\int_{\Gamma_R} e^{w(x)}\left(-\frac{1}{2\pi}\int_{\mathbb R^2_+}F(y)e^{2 w(y)}dy-\frac{1}{2\pi}\int_{\partial \mathbb R^2_+}e^{ w(y)}dl\right)dl\\
    &+\frac{1}{2}\int_{\Gamma_R} e^{ w(x)}\left(-\frac{1}{2\pi}\int_{\mathbb R^2_+}\left(\frac{\langle x+y,x-y\rangle}{|x-y|^2}+\frac{\langle x+\overline y,x-\overline y\rangle}{|\overline x-y|^2}\right)F(y)e^{2w(y)}dy\right.\\
    &\left.-\frac{1}{2\pi}\int_{\partial \mathbb R^2_+}\left(\frac{\langle x+y,x-y\rangle}{|x-y|^2}+\frac{\langle x+\overline y,x-\overline y\rangle}{|\overline x-y|^2}\right)e^{ w(y)}dl\right)dl\\
    = & -\frac{1}{2}d\left(\int_{B_R^+(0)} F(x)e^{2 w(x)}dx+\int_{\Gamma_R}e^wd\l\right)-\frac{1}{4\pi}(I_1+I_2+I_3+I_4).
\end{align*}

To simplify, consider
$$G(x,y):=\frac{\langle x+y,x-y\rangle}{|x-y|^2}+\frac{\langle x+\overline y,x-\overline y\rangle}{|\overline x-y|^2}.$$
Using Fubini's Theorem and changing variables $x$ and $y$, since $G(x,y)=-G(y,x)$ and $G(-x,-y)=G(x,y),$ we get
\begin{align*}
    I_1= & \int_{B_R^+(0)}\int_{\mathbb R^2_+}G(x,y)F(y)e^{2w(y)} F(x)e^{2 w(x)}dydx\\
    = & \int_{B_R^+(0)}\int_{\mathbb R^2_+\backslash B_R^+(0)}G(x,y)F(y)e^{2 w(y)} F(x)e^{2w(x)}dydx\\
    = & \int_{B_{R/2}^+(0)}\int_{\mathbb R^2_+\backslash B_R^+(0)}G(x,y)F(y)e^{2 w(y)} F(x)e^{2w(x)}dydx\\
    & + \int_{B_R^+(0)\backslash B_{R/2}^+(0)}\int_{\mathbb R^2_+\backslash B_{2R}^+(0)}G(x,y)F(y)e^{2 w(y)} F(x)e^{2w(x)}dydx\\
    &+ \int_{B_R^+(0)\backslash B_{R/2}^+(0)}\int_{B_{2R}^+(0)\backslash B_{R}^+(0)}G(x,y)F(y)e^{2 w(y)} F(x)e^{2w(x)}dydx\\
    =& I_{11}+I_{12}+I_{13}.
\end{align*}
If either $|x|\leq R/2$ and $|y|\geq R$ or $R/2\leq|x|\leq R$ and $|y|\geq 2R$, it holds $|G(x,y)|\leq 6$. This implies that
$$|I_{11}| \leq 6\int_{B_{R/2}^+(0)}|F(y)|e^{2w(y)}dy\int_{\mathbb R^2_+\backslash B_R^+(0)}|F(x)|e^{2w(x)}dx$$
and
$$|I_{12}|\leq 6\int_{B_R^+(0)\backslash B_{R/2}^+(0)}|F(y)|e^{2w(y)}dy\int_{\mathbb R^2_+\backslash B_{2R}^+(0)}|F(x)|e^{2w(x)}dx.$$
Since $Fe^{2w}\in L^1(\mathbb R^2_+)$, then $I_{11}$ and $I_{12}$ go to zero as $R\to\infty$.

For $I_{13}$, from Claim 1 we obtain $e^{2w(x)}\leq C|x|^{-2d}$ when $|x|\geq 1$ is sufficiently large. Also for $|x|\geq 1$ we have $|F(x)|\leq C|x|^2$. This implies that
\begin{align}
    \left|\int_{B_{2R}^+(0)\backslash B_{R}^+(0)} G(x,y)F(y)e^{2 w(y)} dy\right|  &\leq  \left|\int_{B_{2R}^+(0)\backslash B_{R}^+(0)}\left(\frac{|x+y|}{|x-y|}+\frac{|x+\overline y|}{|x-\overline y|}\right)|F(y)|e^{2 w(y)} dy\right|\nonumber\\
     &\leq CR^{3-2d}\left|\int_{B_{2R}^+(0)\backslash B_{R}^+(0)}\left(\frac{1}{|x-y|}+\frac{1}{|x-\overline y|}\right) dy\right|\nonumber\\
     &\leq CR^{3-2d}\int_{B_{3R}(0)}\frac{1}{|z|} dz\leq CR^{4-2d}\leq C,\label{eq057}
\end{align}
since $d\geq 2$. This and the fact that $\int_{\mathbb R^2_+}Fe^{2w}dx\in\mathbb R$ implies that $I_{13}$ goes to zero when $R\to\infty$. Therefore, $I_1$ goes to zero when $R\to\infty$.

Similarly, we obtain

\begin{align*}
    I_4 = & \int_{\Gamma_R}\int_{\partial \mathbb R^2_+}G(x,y)e^{w(y)} e^{w(x)}dydx= \int_{\Gamma_R}\int_{\partial \mathbb R^2_+\backslash \Gamma_R}G(x,y)e^{w(y)} e^{w(x)}dydx\\
    = & \int_{\Gamma_{R/2}}\int_{\partial \mathbb R^2_+\backslash \Gamma_R}G(x,y)e^{w(y)} e^{w(x)}dydx +\int_{\Gamma_R\backslash\Gamma_{R/2}}\int_{\partial \mathbb R^2_+\backslash \Gamma_{2R}}G(x,y)e^{w(y)} e^{w(x)}dydx\\
     & +\int_{\Gamma_R\backslash\Gamma_{R/2}}\int_{\Gamma_{2R}\backslash \Gamma_R}G(x,y)e^{w(y)} e^{w(x)}dydx\\
    = &\; I_{41}+I_{42}+I_{43}.
\end{align*}

The estimates for $I_{41}$ and $I_{42}$ are analogous to those for $I_{11}$ and $I_{12}$, but this time we use $e^w\in L^1(\partial\mathbb R^2_+)$. Note that
\begin{align*}
    I_{43} = & \; \int_{R/2}^R\int_R^{2R}G(x,y)\left(e^{w(y)} e^{w(x)}+e^{w(-y)} e^{w(-x)}\right)dydx\\
    & +\int_{R/2}^R\int_R^{2R}G(-x,y)\left(e^{w(-y)} e^{w(x)}+e^{w(y)} e^{w(-x)}\right)dydx.
\end{align*}
Observe that we are making an abuse of notation in the above integrals.

As before, $e^{2w(x)}\leq |x|^{-2\mu}$ for any $\mu<d$ when $|x|\geq 1$ is large enough. This implies that
\begin{align*}
    |I_{43}| \leq & CR^{1-2\mu}\int_{R/2}^R\int_R^{2R}\left(\frac{1}{y-x}+\frac{1}{x+y}\right)dydx\leq CR^{2-2\mu}\ln R.
\end{align*}
For $\mu=3/2$, this implies that $I_{43}$ goes to zero when $R\to\infty$, and therefore $I_4$ goes to zero.

Again, using Fubini's Theorem and changing variables $x$ and $y$, we get
\begin{align*}
    I_2= & \int_{B_R^+(0)}\int_{\partial \mathbb R^2_+}G(x,y)e^{w(y)} F(x)e^{2w(x)}dl dx= -\int_{\partial \mathbb R^2_+}\int_{B_R^+(0)}G(x,y)e^{w(x)} F(y)e^{2w(y)}dydl\\
    = & -\int_{\Gamma_R}\int_{B_R^+(0)}G(x,y)e^{w(x)} F(y)e^{2w(y)}dydl -\int_{\partial \mathbb R^2_+\backslash\Gamma_R}\int_{B_R^+(0)}G(x,y)e^{w(x)} F(y)e^{2w(y)}dydl.
\end{align*}
\begin{align*}
    I_3= & \int_{\Gamma_R}\int_{ \mathbb R^2_+}G(x,y) F(y)e^{2w(y)}e^{ w(x)}dydl\\
    = & \int_{\Gamma_R}\int_{ B_R^+(0)}G(x,y) F(y)e^{2w(y)}e^{w(x)}dydl+ \int_{\Gamma_R}\int_{ \mathbb R^2_+\backslash B_R^+(0)}G(x,y)F(y)e^{2w(y)}e^{w(x)}dydl.
\end{align*}

Thus
\begin{align*}
    I_2+I_3 = & \int_{\Gamma_R}\int_{ \mathbb R^2_+\backslash B_R^+(0)}G(x,y)F(y)e^{2w(y)}e^{w(x)}dydl\\
    &-\int_{\partial \mathbb R^2_+\backslash\Gamma_R}\int_{B_R^+(0)}G(x,y)e^{w(x)} F(y)e^{2w(y)}dydl= II_1-II_2.
\end{align*}

We will employ a technique analogous to the one used to estimate the terms $I_1$ and $I_4$. Write
\begin{align*}
    II_1 = & \int_{\Gamma_{R/2}}\int_{ \mathbb R^2_+\backslash B_R^+(0)}G(x,y)F(y)e^{2w(y)}e^{w(x)}dydl\\
    & + \int_{\Gamma_R\backslash\Gamma_{R/2}}\int_{ \mathbb R^2_+\backslash B_{2R}^+(0)}G(x,y)F(y)e^{2w(y)}e^{w(x)}dydl\\
    & +\int_{\Gamma_R\backslash\Gamma_{R/2}}\int_{ B_{2R}^+(0)\backslash B_R^+(0)}G(x,y)F(y)e^{2w(y)}e^{w(x)}dydl= II_{11}+II_{12}+II_{13}
\end{align*}
and
\begin{align*}
    II_2 = & \int_{\partial \mathbb R^2_+\backslash\Gamma_R}\int_{B_{R/2}^+(0)}G(x,y)e^{w(x)} F(y)e^{2w(y)}dydl\\
    & +\int_{\partial \mathbb R^2_+\backslash\Gamma_{2R}}\int_{B_R^+(0)\backslash B_{R/2}^+(0)}G(x,y)e^{w(x)} F(y)e^{2w(y)}dydl\\
    & + \int_{\Gamma_{2R}\backslash\Gamma_R}\int_{B_R^+(0)\backslash B_{R/2}^+(0)}G(x,y)e^{w(x)} F(y)e^{2w(y)}dydl=II_{21}+II_{22}+II_{23}.
\end{align*}
As we estimated $I_{11}$ and $I_{12}$, we find that the terms $II_{11}$, $II_{12}$, $II_{21}$ and $II_{22}$ go to zero when $R\to\infty$. For $II_{13}$ and $II_{23}$ we use \eqref{eq057} and the fact that $e^w\in L^1(\partial\mathbb R^2_+)$ to see that they go to zero when $R\to \infty$.

Therefore, we conclude that
\begin{equation}\label{eq051}
    RHS\to -\frac{1}{2}d\left(\int_{\mathbb R^2_+} F(x)e^{2 w(x)}dx+\int_{\mathbb R^2_+}e^wdl\right)=-\frac{\pi}{2}d^2.
\end{equation}

If $x=(s,t)$, then we have $\langle x,\nabla w\rangle=s\partial w/\partial s$ for $x\in\Gamma_R$. Thus, using \eqref{eq0000}, we obtain
$$    \int_{\Gamma_{R}}\langle x,\nabla w\rangle\frac{\partial w}{\partial \nu}dl=\int_{\Gamma_{R}}s\frac{\partial w}{\partial s} e^wds =\int_{\Gamma_{R}}\frac{\partial}{\partial s}(se^w)ds-\int_{\Gamma_{R}}e^wds.$$
Thus
\begin{align*}
    LHS = & \int_{B_R^+(0)} F(x)e^{2 w(x)} \langle x,\nabla w\rangle dx+\int_{\Gamma_R}\langle x,\nabla w\rangle e^wdl\\
    = & \frac{1}{2}\int_{B_R^+(0)} F(x) \langle x,\nabla e^{2 w(x)}\rangle dx+\int_{\Gamma_R}\langle x,\nabla w\rangle e^wdl\\
    = & -\frac{1}{2}\int_{B_R^+(0)}e^{2 w(x)}div(F(x)x)dx+\frac{1}{2}R\int_{\partial^+ B_R^+(0)}F(x)e^{2 w(x)}dl+\int_{\Gamma_R}\langle x,\nabla w\rangle e^wdl\\
    = & -\frac{1}{2}\int_{B_R^+(0)}\left(\langle x,\nabla F(x)\rangle+2F(x)\right) e^{2w(x)}dx+\frac{1}{2}R\int_{\partial^+ B_R^+(0)}F(x)e^{2 w(x)}dl\\
    & + \int_{\Gamma_{R}}\frac{\partial}{\partial s}(se^w)ds-\int_{\Gamma_{R}}e^wds.
\end{align*}

Using that $Fe^{2w}\in L^1(\mathbb R^2_+)$ and $e^w\in L^1(\partial\mathbb R^2_+)$, we conclude the existence of a sequence $R_i\to\infty$ such that $R_i\displaystyle\int_{\Gamma_{2,R_i}}Fe^{2w}\rightarrow 0$ and $R_ie^{w(\pm R_i,0)}\to 0$. Therefore, 
\begin{equation}\label{eq052}
    LHS\to -\frac{1}{2}\int_{\mathbb R^2_+}\langle x,\nabla F(x)\rangle e^{2w(x)}dx-\pi d.
\end{equation}
Then, \eqref{eq051} and \eqref{eq052} imply that
$$\frac{1}{\pi}\int_{\mathbb R^2_+}\langle x,\nabla F(x)\rangle e^{2w(x)}dx=d(d-2).$$

Since $F(x)=1+(Ax,x)$ and $A$ is a negative definite matrix, it follows that $\langle x,\nabla F(x)\rangle=2(Ax,x)<0$, which implies that $d(d-2)<0$. However, this contradicts the fact that $d>2$. Therefore, we obtain the result.
\end{proof}

\vspace{1cm}
\noindent\textbf{{Data Availability:}} Data sharing not applicable to this article as no datasets were generated or analyzed during the current manuscript.

\bibliography{references.bib}
\bibliographystyle{acm}
\end{document}